\documentclass[journal,twoside,web]{ieeecolor}

\def\logoname{LOGO-generic-web}
\definecolor{subsectioncolor}{rgb}{0,0.541,0.855}
\usepackage{cite}
\usepackage{units}
\usepackage{amsthm}
\newcounter{thm}
\newtheorem{prob}[thm]{Problem}

\usepackage{pgfplots}
\usepgfplotslibrary{groupplots}

\usepackage{hyperref}
\usepackage{amsmath}
\usepackage{amsfonts}
\usepackage{graphicx} 
\usepackage{amssymb}  
\usepackage{pstricks,pst-plot,psfrag}
\usepackage{units}
\usepackage{float}
\usepackage{dsfont}

\usepackage{lettrine}
\usepackage{import}
\usepackage[acronym]{glossaries}

\usepackage{tikz}
\usetikzlibrary{backgrounds,calc,angles,quotes,arrows.meta,arrows,fit,shapes.geometric,shapes.multipart}
\usetikzlibrary{circuits.ee.IEC}
\usetikzlibrary{shapes.gates.logic.US}

\newacronym{acr:cvt}{CVT}{continuously variable transmission}
\newacronym{acr:CoG}{CoG}{center of gravity}
\newacronym{acr:CoP}{CoP}{center of pressure}
\newacronym{acr:dp}{DP}{dynamic programming}
\newacronym{acr:DoF}{DoF}{degrees of freedom}
\newacronym{acr:ecms}{ECMS}{equivalent consumption minimization strategies}
\newacronym{acr:eltms}{ELTMS}{equivalent lap time minimization strategies}
\newacronym{acr:em}{EM}{electric motor}
\newacronym{acr:es2k}{ES2K}{Energy Storage to Kinetic}
\newacronym{acr:F1}{F1}{Formula 1}
\newacronym{acr:FIA}{FIA}{F\'{e}d\'{e}ration Internationale de l'Automobile}
\newacronym{acr:fgt}{FGT}{fixed-gear transmission}
\newacronym{acr:FD}{FD}{final drive}
\newacronym{acr:ice}{ICE}{internal combustion engine}
\newacronym{acr:k2es}{K2ES}{Kinetic to Energy Storage}
\newacronym{acr:mgu}{MGU}{motor generator unit}
\newacronym{acr:mguh}{MGU-H}{motor generator unit heat}
\newacronym{acr:mguk}{MGU-K}{motor generator unit kinetic}
\newacronym{acr:mpc}{MPC}{model predictive control}
\newacronym[description={energy management strategy}, \glslongpluralkey={energy management strategies},\glsshortpluralkey={EMSs}]{EMS}{EMS}{energy management strategy}%
\newacronym{acr:ODE}{ODE}{ordinary differential equation}
\newacronym{acr:pmp}{PMP}{Pontryagin's Minimum Principle}
\newacronym{acr:pu}{PU}{power unit}
\newacronym[description={powertrain operation}, \glslongpluralkey={powertrain operations},\glsshortpluralkey={POs}]{acr:PO}{PO}{powertrain operation}%
\newacronym{acr:socp}{SOCP}{second-order cone program}
\newacronym{acr:soe}{SoE}{state of energy}

\newcommand{\sN}{\mathbb{N}}

\newcommand{\dtds}{\frac{\textnormal{d}t}{\textnormal{d}s}}
\newcommand{\pushright}[1]{\ifmeasuring@#1\else\omit\hfill$\displaystyle#1$\fi\ignorespaces}
\newcommand{\pushleft}[1]{\ifmeasuring@#1\else\omit$\displaystyle#1$\hfill\fi\ignorespaces}
\makeatother
\newif\ifmargincomments 
\margincommentstrue

\newif\ifextendedversion 
\extendedversiontrue

\ifmargincomments

\else

\fi
\usetikzlibrary{external}
\DeclareMathAlphabet{\mathpzc}{OT1}{pzc}{m}{it}

\newif\ifresponse

\ifresponse

\newcommand{\revTCST}[1]{{\leavevmode\color{blue}#1}}

\else

	\newcommand{\revTCST}[1]{#1}

\fi

\begin{document}

	\title{Optimal Endurance Race Strategies for a Fully Electric Race Car under Thermal Constraints
}

\author{Jorn van Kampen, Thomas Herrmann, Theo Hofman and Mauro Salazar
\thanks{Date of submission \today.
}
\thanks{Jorn van Kampen, Theo Hofman and Mauro Salazar are with the Control Systems Technology section, Department of Mechanical Engineering, Eindhoven University of Technology (TU/e), Eindhoven, 5600 MB, The Netherlands.
		E-mails: {\tt\footnotesize \{j.h.e.v.kampen,t.hofman,m.r.u.salazar\}@tue.nl}} 
\thanks{ Thomas Herrmann is with the Institute of Automotive Technology, Department of Mechanical Engineering, Technical University of Munich (TUM).
		E-mail: {\tt\footnotesize thomas.herrmann@tum.de} }
}


\maketitle
\thispagestyle{plain}
\pagestyle{plain}

\begin{abstract}
This paper presents a bi-level optimization framework to compute the \revTCST{offline} maximum-distance race strategies for a fully electric endurance race car, whilst accounting for the low-level vehicle dynamics and the thermal limitations of the powertrain components.
Thereby, the lower level computes the minimum-stint-time for a given charge time and stint length, whilst the upper level leverages that information to jointly optimize the stint length, charge time and number of pit stops, in order to maximize the driven distance in the course of a fixed-time endurance race.
Specifically, we first extend a convex lap time optimization framework to capture low-level vehicle dynamics and thermal models, and use it to create a map linking the charge time and stint length to the achievable stint time.
Second, we leverage the map to frame the maximum-race-distance problem as a mixed-integer second order conic program that can be efficiently solved in a few seconds to the global optimum with off-the-shelf optimization algorithms.
Finally, we showcase our framework \revTCST{for a simulated} \unit[6]{h} race around the Zandvoort circuit.
Our results show that the optimal race strategy can involve partially charging the battery, and that, compared to the case where the stints are optimized for a fixed number of pit stops, jointly optimizing the stints and number of pit stops can significantly increase the driven distance and hence race performance by several laps.
\end{abstract}

\begin{IEEEkeywords}
	Electric vehicles, Endurance racing, Optimal control, Optimization, Race strategy
\end{IEEEkeywords}

\section{Introduction}
\begin{figure}
	\centering 
	\includegraphics[width=1\linewidth,trim= 0 0 0 0,clip]{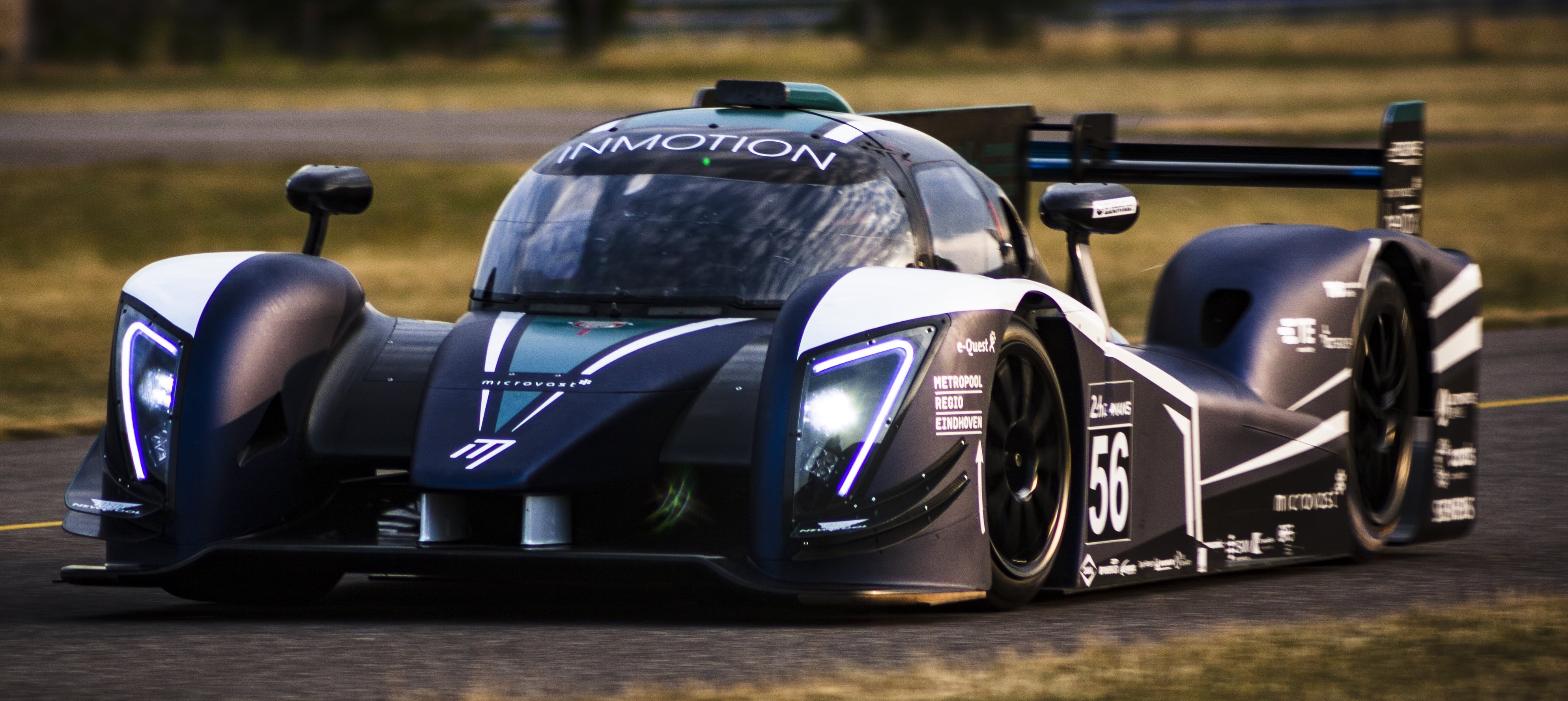}
	\caption{InMotion's fully electric endurance race car~\cite{InMotion}.}
	\label{fig:inmotion-car}
\end{figure}

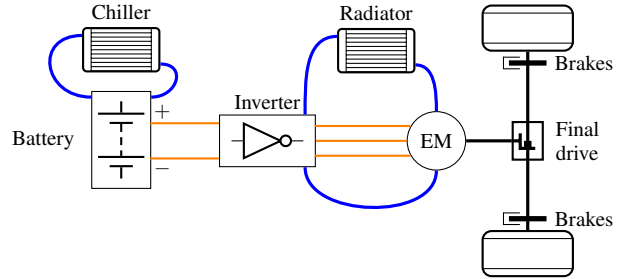
\begin{figure}
	\centering 
	\tikzstyle{simpleNode} = [rectangle, minimum width=1cm, minimum height=1cm,text centered, draw=black, line width=0.25mm]
\tikzstyle{simpleNodeWide} = [rectangle, minimum width=1.5cm, minimum height=1cm,text centered, draw=black, line width=0.25mm]
\tikzstyle{arrow} = [->,-triangle 45,line width=0.25mm,black]
\tikzstyle{lineE} = [-,line width=0.5mm,gray!50,solid]
\tikzstyle{lineM} = [-,line width=0.25mm,black]
\tikzset{
	relative at/.style n args = {3}{
		at = {({$(#1.west)!#2!(#1.east)$} |- {$(#1.south)!#3!(#1.north)$})}
	}
}
\tikzset{BAT/.style={circuit ee IEC}}

%

%
\def\FDwidth{0.4cm} 


%

\begin{tikzpicture}[node distance=2cm]
		
	\tikzstyle{every node}=[font=\footnotesize]; 
	
	\tikzset{trapezium stretches=true};
	\tikzset{three sided/.style={
			draw=none,
			append after command={
				[shorten <= -0.5\pgflinewidth]
				([shift={(-1.5\pgflinewidth,-0.5\pgflinewidth)}]\tikzlastnode.north east)
				edge([shift={( 0.5\pgflinewidth,-0.5\pgflinewidth)}]\tikzlastnode.north west) 
				([shift={( 0.5\pgflinewidth,-0.5\pgflinewidth)}]\tikzlastnode.north west)
				edge([shift={( 0.5\pgflinewidth,+0.5\pgflinewidth)}]\tikzlastnode.south west)            
				([shift={( 0.5\pgflinewidth,+0.5\pgflinewidth)}]\tikzlastnode.south west)
				edge([shift={(-1.0\pgflinewidth,+0.5\pgflinewidth)}]\tikzlastnode.south east)
			}
		}
	}

	\coordinate (Zero) at (0,0);
%

	\node (RR) at (12,3) [rotate=0,draw,thick,rounded corners=1.2mm,minimum width=1.2cm, minimum height=0.6cm] {} node[rectangle,minimum width=1.2cm, minimum height=0.4cm,draw] at (RR) {}; 
	\node (RL) at (12,0) [rotate=0,draw,thick,rounded corners=1.2mm,minimum width=1.2cm, minimum height=0.6cm] {} node[rectangle,minimum width=1.2cm, minimum height=0.4cm,draw] at (RL) {}; 
	\path (RR)--(RL) node (FD) [midway,thick,rectangle,draw, minimum height=0.5cm,minimum width=\FDwidth]  {} node[align=left] at ($(FD.east)+(0.45,0)$) {Final\\drive}; 
	\draw[very thick] (FD.west)--++(0.4pt,0)-- ++(\FDwidth*0.25,0)++(0,0.1cm)--++(0,-0.2cm)--++(\FDwidth*0.5,0) coordinate[midway] (FDgear) ; 
	\draw[line width=0.1cm] (FDgear)--++(0,0.1cm) {};
	\draw[very thick] (FD.south)--(FD.north) {}; 
		
	\draw[very thick] (RR)--(FD.north) {}; 
	\draw[very thick] (RL)--(FD.south) {}; 
	
	\draw[line width=2pt] (RR.south)++(-0.25,-0.15) node (RRb) {} --++ (0.5,0)++(0.5,0) node {Brakes}; 
	\node [three sided,minimum height=1,minimum width=6, inner sep=2,thick] at ($(RRb)+(0.04,0)$) {};
	\draw[line width=2pt] (RL.north)++(-0.25,+0.15) node (RLb) {} --++ (0.5,0)++(0.5,0) node {Brakes}; 
	\node [three sided,minimum height=0.1,minimum width=6,inner sep=2,thick] at ($(RLb)+(0.04,0)$) {};

	
	\node (EM) [circle,radius=10pt,draw] at ($(FD.west) + (-1,0)$) {EM};

	\node (INV) at ($(EM.west) + (-2,0)$) [not gate US,draw,thick] {};
	\draw (INV.input) --++(-0.15,0) node (INVl) {};
	\draw (INV.output) --++(0.15,0) node (INVr) {};
	
	\node[rectangle,fit=(INV) (INVl) (INVr),inner xsep=1pt, inner ysep=5pt,draw] (inverter) {} node[above] at ($(inverter)+(0,0.3)$) {Inverter};
	

	\draw[orange,thick] (EM.150)-|(inverter.18) {};
	\draw[orange,thick] (EM.west)--(inverter.east) {};
	\draw[orange,thick] (EM.210)-|(inverter.-18) {};
	
	\draw[very thick] (EM.east)--(FD.west){};
	

\begin{scope}[circuit ee IEC]
\node at ($(inverter.west)+(-1.3,0.3)$) [battery, point down,thick] (BAT) {} node at ($(BAT.output)+(0,-0.5)$) [battery, point down,thick] (BAT2) {};
\draw[densely dashed,thick] (BAT.output)--(BAT2.input);

\draw[thick] (BAT2.output) --++(0,-0.15) node (BATi) {};
\draw[thick] (BAT.input) --++(0,0.15) node (BATo) {};

\node[rectangle,fit=(BAT) (BATi) (BATo),inner xsep=2pt, inner ysep=0pt,draw] (battery) {} node[left] at ($(battery)+(-0.5,0)$) {Battery};

\draw[orange,thick] (inverter.160) -| (battery.30) ++ (0.15,0.15) node[font=\footnotesize,black] {$+$};
\draw[orange,thick] (inverter.200) -| (battery.-30) ++ (0.15,-0.15) node[font=\footnotesize,black] {$-$};

\end{scope}

\node[rectangle split,minimum width=0.8cm,minimum height=0.5cm,rectangle split every empty part={},draw,rectangle split empty part height=1.2pt, rectangle split parts=10,inner sep=0] (RAD1in) at ($(EM)+(-0.8,1.2)$) {}; 
\node[rectangle,rounded corners=0.5mm,fit=(RAD1in),inner xsep=1mm,inner ysep=0,thick,draw] (RAD1out)  {} node at ($(RAD1out)+(0,0.5)$) {Radiator};

\draw[very thick,blue] (RAD1out.160) to[out=180,in=90,looseness=1] (inverter.35);
\draw[very thick,blue] (inverter.-35) to[in=-90,out=-90,looseness=1] (EM.-90); 
\draw[very thick,blue] (EM.north) to[out=90,in=0,looseness=1]  (RAD1out.-20);

\node[rectangle split,minimum width=0.8cm,minimum height=0.5cm,rectangle split every empty part={},draw,rectangle split empty part height=1.2pt, rectangle split parts=10,inner sep=0] (RAD2in) at ($(battery)+(0,1.2)$) {};
\node[rectangle,rounded corners=0.5mm,fit=(RAD2in),inner xsep=1mm,inner ysep=0,thick,draw] (RAD2out)  {} node at ($(RAD2out)+(0,0.5)$) {Chiller};

\draw[very thick,blue] (battery.125) to[out=180,in=180,looseness=2] (RAD2out.160);
\draw[very thick,blue] (battery.55) to[out=0,in=0,looseness=2] (RAD2out.-20);

\end{tikzpicture}
%
	\caption{Schematic overview of the electric endurance race car powertrain with the electrical connections highlighted in orange and cooling in blue. The battery is actively cooled by a refrigerant system, whereas the electric machine (EM) is cooled by a conventional radiator system.}
	\label{fig:PTlayout}
\end{figure}

\lettrine{T}{he electrification} of race cars has been increasing in popularity over the last years, owing to the advent of hybrid electric Formula 1 cars and Le Mans Hypercars, and battery electric vehicles in Formula E.
In a setting where every millisecond counts, it is of paramount importance to make efficient use of the energy stored on-board via optimized \gls{EMS}, whilst respecting the thermal limits of the powertrain components.
This problem is even more apparent for fully electric race cars, such as the vehicle shown in Fig.~\ref{fig:inmotion-car}, where the available on-board energy is relatively small, compared to conventional race cars. Whereas fuel is the main energy source in conventional race cars, fully electric race cars rely on a less energy-dense battery that is typically integrated within the powertrain according to Fig.~\ref{fig:PTlayout}.

In this context, the necessity of recharging the battery in the course of the race further complicates the problem, requiring race engineers to strike the best trade-off between reducing \revTCST{energy} consumptions and pit stops at the cost of lap time, or driving faster with more pit stops, whilst avoiding damage to the powertrain components by staying within the thermal limits.
This conflict is particularly imminent in endurance racing, where the objective is to maximize the driven distance in a fixed amount of time, which can range up to \unit[24]{h}~\cite{FIA2021}.
In this setting, the car has to be strategically recharged during pit stops in order to maintain a competitive performance and maximize the distance driven.
This calls for algorithms to compute the maximum-distance race strategies that provide the number of pit stops during the race, the number of laps driven per stint (referred to as stint lengths) and the charge time (which is directly correlated with charge energy), whilst accounting for the optimal in-stint strategy in terms of \gls{acr:PO}. Thereby, the \gls{acr:PO} accounts for the thermal management of the powertrain components and for the battery energy management.
Against this backdrop, this paper presents a bi-level optimization framework to compute the maximum-distance race strategies with global optimality guarantees.
  
\subsubsection*{Related Literature}
This work pertains to two main research streams: single-lap optimization of the \glspl{EMS} jointly with the vehicle trajectory or for a given vehicle trajectory, and full-race optimization via simulations.

Several authors have optimized the minimum-lap-time \revTCST{racing line} for a single lap using both direct and indirect optimization methods \revTCST{~\cite{LotEvangelou2013, SedlacekOdenthalEtAl2020, Casanova2000, LimebeerPerantoni2014, ChristWischnewskiEtAl2019, DalBiancoLotEtAl2017,HeilmeierWischnewskiEtAl2020,MassaroLimebeer2021,TremlettLimebeer2016}. As opposed to these qualifying scenarios, racing scenarios are often additionally constrained by energy availability from either the fuel tank or the battery. To approach these racing conditions, some of these single lap studies included a maximum energy consumption per lap~\cite{HerrmannChristEtAl2019,LimebeerPerantoniEtAl2014}}. Similar approaches extend the minimum-lap-time problems to minimum-race-time problems. They consider temperature dynamics, and optimize for multiple consecutive race laps to enable a variable amount of energy consumed per lap, yet formulate the optimization problem in space domain for an a priori known number of laps~\cite{HerrmannPassigatoEtAl2020,LiuFotouhiEtAl2020}. Finally, considering the \revTCST{racing line} to be fixed, multi-lap \glspl{EMS} are optimized, leveraging nonlinear optimization techniques~\cite{HerrmannSauerbeckEtAl2021} or artificial neural networks~\cite{LiuFotouhi2020}. However, these papers lack global optimality guarantees.

Against this backdrop, assuming the \revTCST{racing line} to be available in the form of a maximum speed profile, convex optimization has been successfully leveraged to compute the globally optimal \glspl{EMS} for hybrid and fully electric race vehicles~\cite{EbbesenSalazarEtAl2018,SalazarElbertEtAl2017}, also including gear shift strategies~\cite{DuhrChristodoulouEtAl2020}, different transmission technologies~\cite{BorsboomFahdzyanaEtAl2021} and thermal limitations~\cite{LocatelloKondaEtAl2020}. Yet these methods are either focused on single-lap problems, not capturing pit stops or recharging processes, or they lack global optimality guarantees due to the use of iterative algorithms. 

The second relevant research stream involves race simulations, in which entire races are optimized on a per lap basis~\cite{HeilmeierGrafEtAl2018, WestLimebeer2020}. However, these studies mainly focus on optimal tire strategies by modeling their degradation as a lap time increase and do not capture the \gls{acr:PO} during a lap.

In conclusion, to the best of the authors' knowledge, there are no methods specifically focusing on race strategies in endurance scenarios by which the single-stint operational strategies are jointly optimized \revTCST{together with the high-level stint planning} and accounting for the thermal limits of the powertrain components.

\begin{figure}[!t]
	\centering
	\definecolor{mycolor1}{rgb}{0.00000,0.7500,1.0}%
\definecolor{mycolor2}{rgb}{1.000,0.65000,0.00}%
\definecolor{mycolor3}{rgb}{1.000,0.2000,0.200}%

\definecolor{teal}{rgb}{0.23,0.74,0.7}%
\definecolor{dark-teal}{rgb}{0,0.25,0.36}%
\definecolor{TUeRed}{rgb}{0.78,0.1,0.1}%

\tikzstyle{input} = [rectangle,rounded corners=2mm, minimum width=2cm, minimum height=0.5cm, text centered, draw=black, line width=0.25mm, fill=dark-teal,text=white]
\tikzstyle{output} = [rectangle,rounded corners=2mm, minimum width=2cm, minimum height=0.5cm,text centered, draw=black, line width=0.25mm, fill=teal,text=white]
\tikzstyle{OCP} = [rectangle, minimum width=1cm, minimum height=0.5cm,text centered, draw=black, line width=0.25mm, fill=TUeRed,text=white]
\tikzstyle{LUT} = [rectangle,rounded corners=2mm, minimum width=2cm, minimum height=0.5cm,text centered, draw=black, line width=0.25mm, fill=TUeRed,text=white, on background layer]

\tikzstyle{simpleNodeWide} = [rectangle, minimum width=1.5cm, minimum height=1cm,text centered, draw=black, line width=0.25mm]
\tikzstyle{arrow} = [->,-triangle 45,line width=0.25mm,black]
\tikzstyle{lineE} = [-,line width=0.5mm,gray!50,solid]
\tikzstyle{lineM} = [-,line width=0.25mm,black]
\tikzset{
	relative at/.style n args = {3}{
		at = {({$(#1.west)!#2!(#1.east)$} |- {$(#1.south)!#3!(#1.north)$})}
	}
}

\begin{tikzpicture}[node distance=0.7cm]
	\tikzstyle{every node}=[font=\footnotesize] 
	
	\coordinate (Zero) at (0,0);
	\node (Ict) [input, align=center] at (Zero) {Charge time};
	\node (Isl) [input, align=center] at ($(Ict) + (90:0.75cm)$) {Stint length};
	\node (Irt) [input, align=center] at ($(Isl) + (90:3.3cm)$) {Race time}; 
	\node (HLP) [OCP, align=center] at ($(Irt) + (0:3cm)$) {High-level \\ optimization};
	\node (Osl) [output, align=center] at ($(HLP) + (0:3cm)$) {Stint lengths};
	\node (Ops) [output, align=center, above of=Osl] {Pit stops};
	\node (Oct) [output, align=center, below of=Osl] {Charge times};
	\node (LLP) [OCP, align=center] at ($(Ict) + (0:3cm)$) {Low-level \\ optimization};
	\node (Ovs) [output, align=center]  at ($(LLP) + (0:3cm)$){Vehicle state \\ trajectories};
	
	\draw [-latex,thick] (Irt.east) -- (HLP.west);  
	\draw [-latex,thick] (HLP.east) -- (Ops.west);
	\draw [-latex,thick] (HLP.east) -- (Osl.west);
	\draw [-latex,thick] (HLP.east) -- (Oct.west);
	
	\draw [-latex,thick] (LLP.east) -- (Ovs.west);
	\draw [-latex,thick] (Ict.east) -- (LLP.west);
	\draw [-latex,thick] (Isl.east) -++(0.5,0) |- (LLP.west);
	
	\node[ inner sep=0pt,fill=white] (lookup) at ($(LLP.north)+(90:1.2cm)$) 
	{\resizebox{2.5cm}{!} {\input{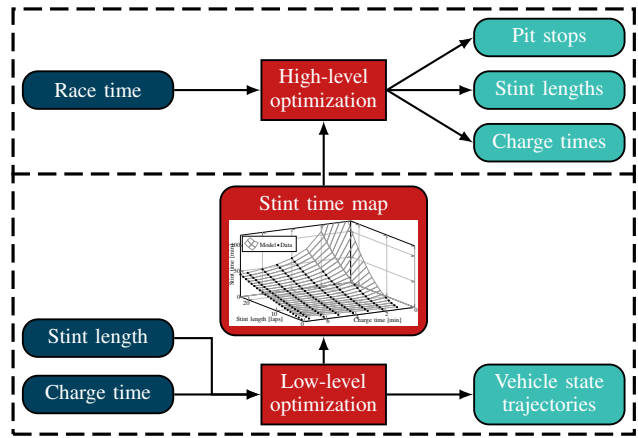}}};
	\node[ align=center, text=white,inner ysep=0] (LUTtext) at ($(lookup.north)+(90:0.2cm)$) {Stint time map}; 
	
	\begin{scope}[on background layer]
			\node [LUT,fit=(lookup) (LUTtext)] (LUTbox) {};
	\end{scope}

	\draw [-latex,thick] (LLP.north) -- (LUTbox.south);
	\draw [-latex,thick] (LUTbox.north) -- (HLP.south);
	
	\node (Box) [rectangle,dash pattern=on 7pt off 3pt, draw=black, very thick, fit = (Irt) (Ops) (Oct) (Ict) (Ovs)] {};
\coordinate (ysplit) at ($(LUTbox.north)!0.5!(Oct.south)$);
\draw [very thick, dash pattern=on 7pt off 3pt] ($(Box.west|-ysplit)$) -- ($(Box.east|-ysplit)$);
	
\end{tikzpicture}                
	\caption{Block diagram overview of the general approach to optimize the race strategy. The low-level framework provides the state trajectories and minimum stint time for several combinations of stint length and charge time. The high-level framework optimizes the race strategy for a given race time, using the stint time data.}
	\label{fig:framework}
\end{figure}

\subsubsection*{Statement of Contributions}
This paper presents a bi-level mixed-integer convex optimization framework to efficiently compute the globally optimal maximum-distance endurance race strategies and the corresponding \gls{acr:PO} in the individual stints, \revTCST{disregarding dynamic events such as overtaking maneuvers, safety car situations and weather changes}.
In order to optimize the strategy over the complete race duration, we decompose an entire endurance race into separate stints and pit stops.
Our low-level algorithm computes the optimal stint time for a given number of laps and level of recharged battery energy. 
Subsequently, we fit the relationship between the stint length, the charged energy, and the achievable stint time as a second-order conic constraint, which we leverage in the high-level algorithm.
Thereby we frame the maximum-distance race problem as a mixed-integer second-order conic program, which jointly optimizes the stint length, the charge time and the number of pit stops.
The resulting problem can be rapidly solved with off-the-shelf numerical solvers with global optimality guarantees.
This procedure is schematically visualized as a block diagram, shown in Fig.~\ref{fig:framework}.

A preliminary version of this paper was accepted for presentation at the 2022 European Control Conference, and selected for direct publication in a special issue of the European Journal of Control~\cite{KampenHerrmannEtAl2022}. Whilst our previous study mainly focused on energy management, we now aim to more accurately represent an endurance racing scenario by accounting for the temperature dynamics of the powertrain components and for the vehicle dynamics. To this end, we include a more accurate battery loss model that captures the dependence on its energy and temperature, and identify a method to model the battery and \gls{acr:em} temperature dynamics in a convex form.
Moreover, we directly include the vehicle dynamics in the low-level control problem by devising a convex framework in the form of a single-track model (shown in Fig.~\ref{fig:bicycle_model}), so that we can better estimate the vehicle capabilities and no longer rely on a pre-computed maximum speed profile. This way, the velocity of the vehicle is directly constrained in the optimization problem via the maximum tire forces. \revTCST{To account for these extensions, we reformulate the high-level control problem by decomposing the race into stints followed by pit stops, instead of stints being preceded by pit stops and we now separately optimize the last stint. }
Finally, we showcase our framework on the Zandvoort circuit for the vehicle shown in Fig.~\ref{fig:inmotion-car}, with the corresponding powertrain architecture shown in Fig.~\ref{fig:PTlayout}.

\subsubsection*{Organization}
The remainder of this paper is structured as follows: Section~\ref{sec:low level} presents the minimum-stint-time control problem, after which Section~\ref{sec:high level} frames the maximum-race-distance control problem. We discuss some of the limitations of our work in Section~\ref{sec:discussion} and showcase our framework for a \unit[6]{h} race in Section~\ref{Results}. Finally, Section~\ref{Conclusion} draws the conclusions and provides an outlook on future research.

\section{Low-level Stint Optimization} \label{sec:low level}
This section illustrates the minimum-stint-time control problem in space domain, since minimizing the stint time given a fixed distance represents the dual problem of maximizing distance within a fixed time. We leverage an existing convex framework~\cite{BroereSalazar2022}, reformulated to a single-track model without steering and side-slip angles, as we do not consider torque-vectoring, and extend it to allow multi-lap optimization, whilst including the \gls{acr:em} and battery temperature dynamics. From the time-optimal control problem, we obtain the minimum stint time for a given stint length and charge time (which is directly related to available battery energy). 

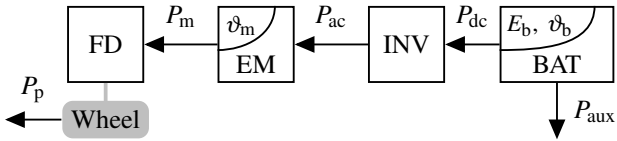
\begin{figure}[!t]
	\centering
	\tikzstyle{simpleNode} = [rectangle, minimum width=1cm, minimum height=1cm,text centered, draw=black, line width=0.25mm]
\tikzstyle{simpleNodeWide} = [rectangle, minimum width=1.5cm, minimum height=1cm,text centered, draw=black, line width=0.25mm]
\tikzstyle{arrow} = [->,-triangle 45,line width=0.25mm,black]
\tikzstyle{lineE} = [-,line width=0.5mm,gray!50,solid]
\tikzstyle{lineM} = [-,line width=0.25mm,black]
\tikzset{
	relative at/.style n args = {3}{
		at = {({$(#1.west)!#2!(#1.east)$} |- {$(#1.south)!#3!(#1.north)$})}
	}
}
\begin{tikzpicture}[node distance=2cm]
\coordinate (Zero) at (0,0);
\node (FD) [simpleNode, align=center] {FD};
\node (EM) [simpleNode, align=center, right of=FD] {};
\node (INV) [simpleNode, align=center, right of=EM] {INV};
\node (BAT) [simpleNodeWide, align=center, right of=INV] {};
\node[relative at={BAT}{0.5}{0.25}] {BAT};
\node[relative at={BAT}{0.35}{0.75}] {\small{$E_\mathrm{b}$, $\vartheta_\mathrm{b}$}};
\node[relative at={EM}{0.5}{0.25}] {EM};
\node[relative at={EM}{0.30}{0.75}] {\small{$\vartheta_\mathrm{m}$}};
\node (Wheel) [simpleNode, align=center, below of=FD, yshift=1cm, minimum height=0.5cm, draw=none, fill=gray!50, rounded corners] {Wheel};
\draw [lineM] ($(BAT.west)+(0,-0.075)$) to [out=0,in=270] ($(BAT.north)+(0.50,0)$);
\draw [lineM] ($(EM.west)+(0,-0.075)$) to [out=0,in=270] ($(EM.north)+(0.25,0)$);
\draw [lineE] (FD.south) -- node[right] {} ($(Wheel.north)-(0,0.1)$);
\draw [arrow] (Wheel.west) -- node[above, yshift=0.1cm] {$P_\mathrm{p}$} ($(Wheel.west) + (-0.75,0)$);
\draw [arrow] (EM.west) -- node[above, yshift=0.1cm] {$P_\mathrm{m}$} (FD.east);
\draw [arrow] (INV.west) -- node[above, yshift=0.1cm] {$P_\mathrm{ac}$} (EM.east);
\draw [arrow] (BAT.west) -- node[above, yshift=0.1cm] {$P_\mathrm{dc}$} (INV.east);
\draw [arrow] (BAT.south) -- node[right, xshift=0.1cm] {$P_\mathrm{aux}$} ($(BAT.south) - (0,0.75)$);
\end{tikzpicture}                
	\caption{Schematic layout of the electric race car powertrain topology consisting of a battery (BAT), inverter (INV), electric machine (EM) and final drive (FD). The arrows indicate positive power flows of the auxiliary power $P_\mathrm{aux}$, the electrical inverter power $P_\mathrm{dc}$, the electrical \gls{acr:em} input power $P_\mathrm{ac}$ and mechanical output power $P_\mathrm{m}$, and the propulsion power $P_\mathrm{p}$. \revTCST{The vehicle is rear-wheel driven, with mechanical brakes on each wheel.}}
	\label{fig:topology}
\end{figure}

Fig.~\ref{fig:topology} shows a schematic representation of the powertrain topology of the electric race car presented in Fig.~\ref{fig:PTlayout}. The EM propels both of the rear wheels through a fixed \gls{acr:FD}, while receiving energy from the battery pack via the inverter. As with most electric vehicles, the EM can also operate as a generator, thus we account for a bi-directional energy flow between the battery and the wheels. In addition, we consider auxiliary components that are powered from the main battery, \revTCST{such as pumps and compressors for the cooling system,} as a uni-directional energy flow.  

\begin{figure}[!t]
	\centering
	\begin{tikzpicture}
	
	\pgfmathsetmacro{\y}{4.9}
	\pgfmathsetmacro{\x}{2}
	\pgfmathsetmacro{\Angle}{0}
	\pgfmathsetmacro{\AngleR}{-20}
	\pgfmathsetmacro{\AngleF}{-35}
	\pgfmathsetmacro{\AngleB}{-25}
	
	\pgfmathsetmacro{\AngleDelta}{-20}
	\pgfmathsetmacro{\COMradius}{0.1}
	
	\tikzstyle{every node}=[font=\footnotesize] 
	
	\definecolor{Fcol}{rgb}{1,0.00,0.00} 
	\definecolor{vcol}{rgb}{0,0.50,0.00} 
	
	\coordinate (Origin) at (0,0);
	\coordinate (yaxis) at (0,\x+1.5);
	\coordinate (xaxis) at (\y+1.5,0);
	\coordinate (Fyf) at (\y,\x); 
	\coordinate (Fyr) at (0,\x); 
	\coordinate (AngleStart) at (\y,\x); 
	\draw [dashed,thick] (AngleStart)--++(\Angle:1.7cm) coordinate (AngleEnd);
	
	\node at (Fyr) [rotate=\Angle,draw,thick,rounded corners=2mm,minimum width=1.2cm, minimum height=0.6cm] {};
	\draw [vcol,-latex,thick] (Fyr)--++(\Angle+\AngleR:1.3cm) coordinate (RedArrowOne) node[below] {$v_\mathrm{R}$};
	\draw pic["$\alpha_\mathrm{R}$", draw=black, text=black, -latex, angle eccentricity=1.45, angle radius=0.8cm]
	{angle=RedArrowOne--Fyr--Fyf};
	\draw [-latex,thick,Fcol] (Fyr)--++(\Angle-90:0.7cm) node [rotate=\Angle,left] {$F_\mathrm{y,R}$};           
	\draw [thick] (Fyr)--(Fyf); 
	
	\node at (Fyf) [rotate=\Angle+\AngleDelta,draw,thick,rounded corners=2mm,minimum width=1.2cm, minimum height=0.5cm] {};         
	\draw [vcol,-latex,thick] (Fyf)--++(\Angle+\AngleF:1.3cm) coordinate (RedArrowTwo) node at ++(-0.5,0) {$v_\mathrm{F}$}; 
	\draw [dashed,thick] (Fyf)--++(\Angle+\AngleDelta:1.3cm) coordinate (DeltaAngleEnd);
	\draw pic["$\delta_\mathrm{F}$", draw=black, text=black, -latex, angle eccentricity=1.35, angle radius=0.8cm]
	{angle=DeltaAngleEnd--Fyf--AngleEnd};
	\draw pic["$\alpha_\mathrm{F}$", draw=black, text=black, -latex, angle eccentricity=1.45, angle radius=1cm]
	{angle=RedArrowTwo--Fyf--DeltaAngleEnd};
	\draw [-latex,thick,Fcol] (Fyf)--++(\AngleDelta-90:0.7cm) node [rotate=\Angle,left] {$F_\mathrm{y,F}$};       
	
	\coordinate (COM) at (0.48*\y,\x); 
	\begin{scope}[rotate=\Angle]
		\fill [radius=\COMradius] (COM) -- ++(\COMradius,0) arc [start angle=0,end angle=90] -- ++(0,-2*\COMradius) arc [start angle=270, end angle=180];
		\draw [thick,radius=\COMradius] (COM) circle node[above] {$\mathrm{CoG}$};
	\end{scope}

	\draw [-latex,thick,vcol] (COM)--++(\AngleB:1cm) coordinate (v) node [below,rotate=\Angle] {$v$};
	\draw pic["$\beta$", right, draw=black, text=black, -latex, angle eccentricity=1.35, angle radius=0.65cm]
	{angle=v--COM--Fyf};
	
	\coordinate (LrLabel) at ($ (Fyr) +  (\Angle+90:0.5cm) $);
	\coordinate (COMLabel) at ($ (COM) +  (\Angle+90:0.5cm) $);
	\coordinate (LfLabel) at ($ (Fyf) +  (\Angle+90:0.5cm) $);
%
\end{tikzpicture}  
	\input{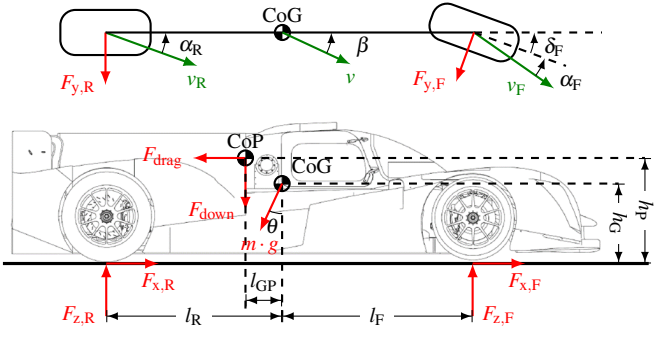}               
	\caption{Top and side view of the single-track vehicle model. The forces acting on the vehicle are shown in red and the velocity vectors are shown in green. The center of gravity (CoG) is shown in both figures, whereas the center of pressure is only shown in the side view. The arrows indicate positive directions, angles or forces.}
	\label{fig:bicycle_model}
\end{figure}

In reality, the driver controls the EM torque and the brake force through the accelerator and brake pedal, respectively, and as such we define the mechanical EM power $P_{\mathrm{m}}$ and mechanical brake power $P_\mathrm{brake}$ as the input variables. As state variables, we choose the battery energy $E_{\mathrm{b}}$, battery temperature $\vartheta_\mathrm{b}$, \gls{acr:em} temperature $\vartheta_\mathrm{m}$ and the kinetic energy of the vehicle $E_{\mathrm{kin}}$. The remaining energy flows between the powertrain components are the propulsion power $P_{\mathrm{p}}$, electrical EM power $P_{\mathrm{ac}}$, electrical inverter power $P_{\mathrm{dc}}$ and auxiliary supply $P_{\mathrm{aux}}$. Since we formulate the control problem in space domain, we ultimately define the model in terms of forces rather than power. Hence, we divide power by the vehicle velocity, since the space-derivative of energy is defined with respect to the vehicle.

\subsection{Objective and Longitudinal Dynamics}
In racing, the objective is to minimize the lap times over the entire race. Since we only consider a stint in the low-level control problem, the objective is to minimize the stint time $t_{\mathrm{stint}}$, which is defined as
\par\nobreak\vspace{-5pt}
\begingroup
\allowdisplaybreaks
\begin{small}
	\begin{equation}
		\min t_{\mathrm{stint}} = \min \int_{0}^{S_{\mathrm{stint}}} \frac{\mathrm{d}t}{\mathrm{d}s}(s) \ \mathrm{d}s,
		\label{eq:low_obj}
	\end{equation}
\end{small}%
\endgroup
where $S_{\mathrm{stint}}$ is the stint length in terms of distance and $\frac{\mathrm{d}t}{\mathrm{d}s}(s)$ is the lethargy, which is the inverse of the vehicle velocity \revTCST{$v(s)\geq v_\mathrm{min}$, with $v_\mathrm{min}$ being a positive velocity close to standstill}. To implement the lethargy as a convex constraint, we define
\par\nobreak\vspace{-5pt}
\begingroup
\allowdisplaybreaks
\begin{small}
	\begin{equation}
		\frac{\mathrm{d}t}{\mathrm{d}s}(s) \geq \frac{1}{v(s)},
		\label{eq:lethargy}
	\end{equation}
\end{small}%
\endgroup
which is a convex relaxation that holds with equality in case of an optimal solution~\cite{EbbesenSalazarEtAl2018}.

In this study, we limit ourselves to non-torque-vectoring powertrain topologies and thereby only include the longitudinal vehicle dynamics with the use of a single-track model that captures the front and rear axle individually, as shown in Fig.~\ref{fig:bicycle_model}.\revTCST{ The single-track model was found to provide the best trade-off between model complexity and accuracy for race strategy optimization purposes, and has been frequently applied in minimum lap time studies~\cite{MassaroLimebeer2021}}. The longitudinal force balance is given by
\par\nobreak\vspace{-5pt}
\begingroup
\allowdisplaybreaks
\begin{small}
	\begin{equation}
		\frac{\mathrm{d}}{\mathrm{d}s}E_{\mathrm{kin}}(s) = F_{\mathrm{x,F}}(s)+F_{\mathrm{x,R}}(s)		-F_{\mathrm{drag}}(s)-m\cdot g \cdot \sin(\theta(s)),
		\label{eq:longitudinal}
	\end{equation}
\end{small}%
\endgroup
where $F_{\mathrm{x},i}(s)$ is the longitudinal force per axle with $i \in [\mathrm{F, R}]$ denoting the front and rear axle, respectively, $F_{\mathrm{drag}}(s)$ is the aerodynamic drag force, $m$ is the total mass of the vehicle, $g$ is the gravitational constant and $\theta(s)$ is the inclination of the track along the racing line. \revTCST{Hereby we assume a track with a fixed racing line, and known altitude and banking profiles.} The aerodynamic drag force is given by
\par\nobreak\vspace{-5pt}
\begingroup
\allowdisplaybreaks
\begin{small}
	\begin{equation}
		F_{\mathrm{drag}}(s)=\frac{c_\mathrm{d} \cdot A_{\mathrm{f}} \cdot \rho}{m} \cdot E_{\mathrm{kin}}(s),
		\label{eq:dragforce}
	\end{equation}
\end{small}%
\endgroup
where $c_\mathrm{d}$ is the drag coefficient, $A_\mathrm{f}$ is the frontal area of the vehicle and $\rho$ is the air density
The longitudinal axle forces are defined as
\par\nobreak\vspace{-5pt}
\begingroup
\allowdisplaybreaks
\begin{small}
	\begin{equation}
		F_{\mathrm{x},i}(s)=F_{\mathrm{p},i}(s)-c_{\mathrm{r}}\cdot F_{\mathrm{z},i}(s)-F_{\mathrm{brake},i}(s),
		\label{eq:Fx}
	\end{equation}
\end{small}%
\endgroup
where $F_{\mathrm{p},i}(s)$ is the propulsion force per axle, $c_{\mathrm{r}}$ is the rolling resistance coefficient, $F_{\mathrm{z},i}(s)$ represents the vertical axle force and $F_{\mathrm{brake},i}(s)$ is the force from the mechanical brakes per axle. For non-driven wheels we set $F_{\mathrm{p},i}(s)=0$, whereas for driven wheels we write~\eqref{eq:Fx} as two inequality constraints to capture the final drive losses through
\par\nobreak\vspace{-5pt}
\begingroup
\allowdisplaybreaks
\begin{small}
	\begin{equation} \label{eq:Fx_trac}
		F_{\mathrm{x},i}(s) \leq F_{\mathrm{m}}(s)\cdot \eta_{\mathrm{fd}}-c_{\mathrm{r}}\cdot F_{\mathrm{z},i}(s)-F_{\mathrm{brake},i}(s),
	\end{equation} 
\end{small}%
\endgroup
\par\nobreak\vspace{-5pt}
\begingroup
\allowdisplaybreaks
\begin{small}
	\begin{equation} \label{eq:Fx_reg}
		F_{\mathrm{x},i}(s) \leq F_{\mathrm{m}}(s)\cdot \frac{1}{\eta_{\mathrm{fd}}}-c_{\mathrm{r}}\cdot F_{\mathrm{z},i}(s)-F_{\mathrm{brake},i}(s),
	\end{equation}
\end{small}%
\endgroup
where $F_{\mathrm{m}}(s)$ is the mechanical output force from the EM and $\eta_{\mathrm{fd}}$ is the efficiency of the final drive, assumed constant. Due to the objective~\eqref{eq:low_obj}, in case of traction, \eqref{eq:Fx_trac} will hold with equality, whilst in the case of regenerative braking, \eqref{eq:Fx_reg} will hold with equality, thus capturing the bi-directional power flow.

The lateral force balance is defined as 
\par\nobreak\vspace{-5pt}
\begingroup
\allowdisplaybreaks
\begin{small}
	\begin{equation}
		2 E_{\mathrm{kin}}(s)\cdot \zeta(s) = F_{\mathrm{y,F}}(s) + F_{\mathrm{y,R}}(s) \revTCST{ + m\cdot g\cdot \sin(\phi(s))},
		\label{eq:lateral}
	\end{equation}
\end{small}%
\endgroup
where $\kappa(s)$ is the pre-computed curvature of the optimal \revTCST{racing line} around the track, assumed to be fixed, $F_{\mathrm{y},i}(s)$ represents the lateral force per axle \revTCST{and $\phi(s)$ is the banking of the track along the racing line}. 

The vertical force balance consists of the static load and the aerodynamic downforce (cf. Fig.~\ref{fig:bicycle_model}) as
\par\nobreak\vspace{-5pt}
\begingroup
\allowdisplaybreaks
\begin{small}
	\begin{equation}
		F_{\mathrm{z,F}}(s) + F_{\mathrm{z,R}}(s) = m\cdot g\cdot \revTCST{\cos(\theta(s))\cdot\cos(\phi(s))} + F_{\mathrm{down}}(s),
		\label{eq:vertical}
	\end{equation}
\end{small}%
\endgroup
where $F_{\mathrm{z},i}(s)$ is the vertical force per axle and $ F_{\mathrm{down}}(s)$ is the aerodynamic downforce given by
\par\nobreak\vspace{-5pt}
\begingroup
\allowdisplaybreaks
\begin{small}
	\begin{equation}
		F_{\mathrm{down}}(s)=\frac{c_\mathrm{l} \cdot A_{\mathrm{f}} \cdot \rho}{m} \cdot E_{\mathrm{kin}}(s),
		\label{eq:downforce}
	\end{equation}
\end{small}%
\endgroup
where $c_\mathrm{l}$ is the lift coefficient.
We consider steady-state cornering only, thereby assuming a yaw moment equilibrium given by
\par\nobreak\vspace{-5pt}
\begingroup
\allowdisplaybreaks
\begin{small}
	\begin{equation}
		F_{\mathrm{y,F}}(s)\cdot l_{\mathrm{F}} = F_{\mathrm{y,R}}(s)\cdot l_{\mathrm{R}},
		\label{eq:yaw-moment}
	\end{equation}
\end{small}%
\endgroup
where $l_{i}$ represents the horizontal distance from the respective axle to the \gls{acr:CoG}. 

The longitudinal load transfer is determined through the pitch moment equilibrium \revTCST{about the orthogonal projection of the \gls{acr:CoG} on the road plane}, which is defined by
\par\nobreak\vspace{-5pt}
\begingroup
\allowdisplaybreaks
\begin{small}
	\begin{multline}
		\frac{\mathrm{d}}{\mathrm{d}s}E_{\mathrm{kin}}(s)\cdot h_{\mathrm{G}} =
		F_{\mathrm{z,R}}(s)\cdot l_{\mathrm{R}} - F_{\mathrm{z,F}}(s)\cdot l_{\mathrm{F}} -m\cdot g \cdot \sin(\theta(s))\cdot h_{\mathrm{G}} \\
		- F_{\mathrm{drag}}(s)\cdot h_{\mathrm{P}} - F_{\mathrm{down}}(s)\cdot l_{\mathrm{GP}},
		\label{eq:pitch-moment}
	\end{multline}
\end{small}%
\endgroup
where $h_{\mathrm{G}}$ is the height of the \gls{acr:CoG} with respect to the ground, $h_{\mathrm{P}}$ is the height of the \gls{acr:CoP} with respect to the ground and $l_{\mathrm{GP}}$ is the horizontal distance from the \gls{acr:CoG} to the \gls{acr:CoP}. 

The longitudinal and lateral forces are bounded by their respective friction circles per axle, which are defined by the convex set written as\revTCST{
\par\nobreak\vspace{-5pt}
\begingroup
\allowdisplaybreaks
\begin{small}
	\begin{equation}
		\left(\frac{F_{\mathrm{x},i}(s)}{\mu_{\mathrm{x},i}}\right)^2 +\left( \frac{F_{\mathrm{y},i}(s)}{\mu_{\mathrm{y},i}}\right)^2 \leq F_{\mathrm{z},i}(s)^2,
		\label{eq:friction_circle}
	\end{equation}
\end{small}%
\endgroup
where $\mu_{\mathrm{x},i}$ and $\mu_{\mathrm{y},i}$, both assumed constant, represent the longitudinal and lateral tire friction coefficients, respectively}. Although the constraint function is not convex, it specifies a convex set, which is shown in the Appendix.

The majority of racing vehicles are equipped with mechanical brakes that provide a fixed brake force ratio between the front and rear wheels. Therefore, we define a relation between the front and rear brake force as
\par\nobreak\vspace{-5pt}
\begingroup
\allowdisplaybreaks
\begin{small}
	\begin{equation}
		F_{\mathrm{brake,R}}(s)\cdot \delta_{\mathrm{brake}} = F_{\mathrm{brake,F}}(s)\cdot (1 - \delta_{\mathrm{brake}}),
		\label{eq:brake_distribution}
	\end{equation}
\end{small}%
\endgroup
where $\delta_{\mathrm{brake}}$ represents the brake balance with respect to the front, which is assumed to remain constant during the race. \revTCST{Since we are interested in finding the optimal powertrain operation, we decouple the mechanical brakes and the regenerative braking from the \gls{acr:em}. In practice, this can be regarded as a separate driver input for regenerative braking, but the exact implementation is beyond the scope of this work.}

The relation between the kinetic energy and velocity of the vehicle is defined by a convex relaxation as
\par\nobreak\vspace{-5pt}
\begingroup
\allowdisplaybreaks
\begin{small}
	\begin{equation} \label{eq:E_kin}
		E_{\mathrm{kin}}(s) \geq \frac{1}{2} \cdot m\cdot v^2(s).
	\end{equation}
\end{small}%
\endgroup
In contrast to single-lap scenarios, a stint is represented by the vehicle starting and stopping at the pit box with a certain number of flying laps in between. However, since we are working in space domain, the lethargy would diverge to infinity for zero velocity. To solve this issue, we define a minimal velocity $v_{\mathrm{min}}$ close to standstill and enforce this value to the initial and final velocity with
\par\nobreak\vspace{-5pt}
\begingroup
\allowdisplaybreaks
\begin{small}
	\begin{equation} \label{eq:E0_ES}
		E_{\mathrm{kin}}(0) = E_{\mathrm{kin}}(S_{\mathrm{stint}}) = \frac{1}{2}\cdot m\cdot v_{\mathrm{min}}^2.
	\end{equation}
\end{small}%
\endgroup
Finally, the vehicle should adhere to a strict speed limit, of which the exact value is track-dependent, when driving through the pit lane. Therefore, we define an upper bound $v_{\mathrm{pit,max}}$ on the vehicle velocity when the vehicle is exiting or entering the pit as
\par\nobreak\vspace{-5pt}
\begingroup
\allowdisplaybreaks
\begin{small}
	\begin{equation} \label{eq:Ek_pit}
		E_{\mathrm{kin}}(s) \leq \frac{1}{2}\cdot m \cdot v_{\mathrm{pit,max}}^2 \quad \forall s \in \mathpzc{S}_{\mathrm{pit}},
	\end{equation}
\end{small}%
\endgroup
where $\mathpzc{S}_{\mathrm{pit}}$ is the set of distance-based positions that are part of the pit lane.

\subsection{Electric Machine}\label{sec:EM}
This section derives a convex representation of the operating limits and power losses of the \gls{acr:em}. Moreover, we derive the thermal dynamics of the \gls{acr:em} and compare the model against real-world test data.  

In general, we can distinguish between a maximum torque and maximum power operating region for an \gls{acr:em}. Translating this to constraints in space domain results in a lower and upper bound on the mechanical output force of the EM for the maximum torque region as
\par\nobreak\vspace{-5pt}
\begingroup
\allowdisplaybreaks
\begin{small}
	\begin{equation} \label{eq:Tmax}
		F_{\mathrm{m}}(s) \in \left[-\frac{T_{\mathrm{m,max}}\cdot \gamma_{\mathrm{fd}}}{r_{\mathrm{w}}},\frac{T_{\mathrm{m,max}}\cdot \gamma_{\mathrm{fd}}}{r_{\mathrm{w}}}\right],
	\end{equation}
\end{small}%
\endgroup
where $T_{\mathrm{m,max}}$ is the maximum torque the \gls{acr:em} can deliver, $\gamma_{\mathrm{fd}}$ is the final drive ratio and $r_{\mathrm{w}}$ is the radius of the rear wheels. Note that we include the final drive ratio, as we define the space-derivatives with respect to the vehicle reference frame. Similarly, the mechanical output force of the \gls{acr:em} within the maximum power region is bounded as
\par\nobreak\vspace{-5pt}
\begingroup
\allowdisplaybreaks
\begin{small}
	\begin{equation} 
		F_{\mathrm{m}}(s) \in \left[-P_{\mathrm{m,max}}\cdot \frac{\mathrm{d}t}{\mathrm{d}s}(s), P_{\mathrm{m,max}}\cdot \frac{\mathrm{d}t}{\mathrm{d}s}(s) \right],
		\label{eq:P_em,max}
	\end{equation}
\end{small}%
\endgroup
where $P_{\mathrm{m,max}}$ is the maximum power the \gls{acr:em} can deliver. 

We model the \gls{acr:em} \emph{force} losses $F_{\mathrm{m,l}}(s)$ rather than the power losses as a function of the vehicle velocity and force of the EM. \revTCST{We calculate the force loss data points, denoted with $(\hat{.})$, from the \gls{acr:em} power losses as
\par\nobreak\vspace{-5pt}
\begingroup
\allowdisplaybreaks
\begin{small}
	\begin{equation} 
		\hat{F}_{\mathrm{m,l}} = \frac{\hat{P}_\mathrm{m,l}\cdot \gamma_{\mathrm{fd}}}{\hat{\omega}_\mathrm{m}\cdot r_\mathrm{w}}.
		\label{eq:P_m,loss}
	\end{equation}
\end{small}%
\endgroup}
 In general, an \gls{acr:em} efficiency map shows large losses at low rotational velocities. Hence, we want to include a term in our losses fit that is inversely proportional to the vehicle velocity. To ensure convexity, we model the \gls{acr:em} losses as
\par
\begin{figure}
	\vspace{4mm}
	\centering
	\input{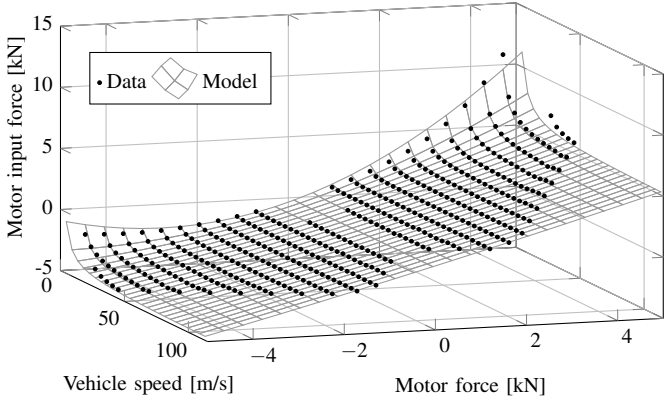}
	\caption{A speed- and torque-dependent model of the EM. The normalized root-mean-square error (RMSE) of the model is 1.49\% w.r.t.\ the maximum motor input force $F_\mathrm{ac}$.  }
	\label{fig:emlossfit}
\end{figure}

\nobreak\vspace{-5pt}
\begingroup
\allowdisplaybreaks
\begin{small}
	\begin{equation}
		F_{\mathrm{m,l}}(s) = x_{\mathrm{m,l}}^{\top}(s) Q_{\mathrm{m,l}}x_{\mathrm{m,l}}(s),
	\end{equation}
\end{small}%
\endgroup
where $x_{\mathrm{m,l}}(s) = \left[\frac{1}{\sqrt{v(s)}} \  \sqrt{v(s)} \  \frac{F_{\mathrm{m}}(s)}{\sqrt{v(s)}}\right]^\top$ and $Q_{\mathrm{m,l}} \in \mathbb{S}_+^3$ is a symmetric positive semi-definite matrix of coefficients, whose values are determined through semi-definite programming. Fig.~\ref{fig:emlossfit} shows the \gls{acr:em} input force as a function of the \gls{acr:em} output force and vehicle speed for the convex model and for the reference data.
To implement the losses in a convex manner, we take the relation of the electrical \gls{acr:em} input force $F_{\mathrm{ac}}(s)$ to the mechanical output force as
\par\nobreak\vspace{-5pt}
\begingroup
\allowdisplaybreaks
\begin{small}
	\begin{equation}
		F_{\mathrm{ac}}(s) = F_{\mathrm{m}}(s) + F_{\mathrm{m,l}}(s),
	\end{equation}
\end{small}%
\endgroup
substitute the loss model, relax it and rewrite to a relaxation describing a convex set as
\par\nobreak\vspace{-5pt}
\begingroup
\allowdisplaybreaks
\begin{small}
	\begin{equation} \label{eq:eff_EM}
		(F_{\mathrm{ac}}(s)-F_{\mathrm{m}}(s))\cdot v(s) \geq y_{\mathrm{m,l}}^{\top}(s) Q_{\mathrm{m,l}}y_{\mathrm{m,l}}(s),
	\end{equation}
\end{small}%
\endgroup
where $y_{\mathrm{m,l}}(s) = \left[1 \ v(s) \  F_{\mathrm{m}}(s)\right]^\top$. The convexity of this constraint is shown in the Appendix by writing it as a second-order conic constraint.
\par

For the cooling circuit of the \gls{acr:em}, we consider a conventional setup using liquid cooling and radiators, as commonly applied in motorsport. The losses are assumed to be converted to heat, thereby changing the \gls{acr:em} temperature according to the first-order temperature \gls{acr:ODE} given by
\par\nobreak\vspace{-5pt}
\begingroup
\allowdisplaybreaks
\begin{small}
	\begin{equation}
		C_{\mathrm{m}}\cdot \frac{\mathrm{d}\vartheta_{\mathrm{m}}}{\mathrm{d}t}(s) = P_{\mathrm{m,l}}(s) - P_{\mathrm{m,c}}(s),
		\label{eq:EM_thermal_time}
	\end{equation}
\end{small}%
\endgroup
where $C_{\mathrm{m}}$ is the total lumped thermal capacity of the \gls{acr:em}, $\vartheta_{\mathrm{m}}(s)$ is the temperature of the \gls{acr:em} (defined at the windings), $P_{\mathrm{m,l}}(s)$ are the \gls{acr:em} power losses and $P_{\mathrm{m,c}}(s) \geq 0$ represents the power outflow to the cooling liquid as
\par\nobreak\vspace{-5pt}
\begingroup
\allowdisplaybreaks
\begin{small}
	\begin{equation}
		P_{\mathrm{m,c}}(s) =\revTCST{ ({\vartheta_{\mathrm{m}}(s)-\vartheta_{\mathrm{m,c}}})\cdot {\lambda_{\mathrm{m}}}},
		\label{eq:EM_cooling_time}
	\end{equation}
\end{small}%
\endgroup
where $\vartheta_{\mathrm{m,c}}$ represents the temperature of the cooling liquid and \revTCST{$\lambda_{\mathrm{m}}$ is the thermal conductivity} between the \gls{acr:em} and the cooling liquid, where we assume both parameters to be constant. In reality, the coolant temperature may change depending on the performance of the cooling system. Yet in high-level motorsports, this change in coolant temperature is relatively small compared to the change in \gls{acr:em} temperature and therefore we deem the modeling of the underlying cooling system beyond the scope of this paper.
Rewriting~\eqref{eq:EM_thermal_time} to space domain results in
\par\nobreak\vspace{-5pt}
\begingroup
\allowdisplaybreaks
\begin{small}
	\begin{equation}
		C_{\mathrm{m}}\cdot \frac{\mathrm{d}\vartheta_{\mathrm{m}}}{\mathrm{d}s}(s) = F_{\mathrm{m,l}}(s) - F_{\mathrm{m,c}}(s),
		\label{eq:EM_thermal_space}
	\end{equation}
\end{small}%
\endgroup
where $ F_{\mathrm{m,c}}(s) \geq 0$ is the force-equivalence of the \gls{acr:em} cooling power. This \gls{acr:em} cooling force is ultimately obtained by rewriting~\eqref{eq:EM_cooling_time} to forces using a convex relaxation and linear equality constraint. To obtain a convex representation, we approximate the \gls{acr:em} temperature as
\par
\begin{figure}[!tb]
	\centering
%
%
\definecolor{mycolor1}{rgb}{0.500000,0.50,0.500}%
\definecolor{mycolor2}{rgb}{0.000,0.00,0.00}%
\begin{tikzpicture}
\footnotesize
\begin{groupplot}[%
	group style={
		group size=1 by 2,
		vertical sep=7pt},
	width=0.86\columnwidth,
	at={(0,0)},
	scale only axis,
	xmin=30,
	xmax=170,
	y label style={yshift=-0.4cm},
	axis background/.style={fill=white},
	xmajorgrids,
	ymajorgrids,
	xlabel={Real temperature [$^\circ$C]},
	x label style={yshift=0.2cm}
	]
\nextgroupplot[%
height=1.7cm,
xticklabels={,,},
ylabel style={align=center},
ylabel={Temperature fit [$^\circ$C]},
ymin=20,
ymax=180,
legend style={at={(0.05,0.9)}, anchor=north west, legend cell align=left, align=left},
legend entries={Data, Model},
legend columns=1,
]
\addplot [color=mycolor1, dashed, line width=0.75pt]
  table[row sep=crcr]{%
30	30\\
31	31\\
32	32\\
33	33\\
34	34\\
35	35\\
36	36\\
37	37\\
38	38\\
39	39\\
40	40\\
41	41\\
42	42\\
43	43\\
44	44\\
45	45\\
46	46\\
47	47\\
48	48\\
49	49\\
50	50\\
51	51\\
52	52\\
53	53\\
54	54\\
55	55\\
56	56\\
57	57\\
58	58\\
59	59\\
60	60\\
61	61\\
62	62\\
63	63\\
64	64\\
65	65\\
66	66\\
67	67\\
68	68\\
69	69\\
70	70\\
71	71\\
72	72\\
73	73\\
74	74\\
75	75\\
76	76\\
77	77\\
78	78\\
79	79\\
80	80\\
81	81\\
82	82\\
83	83\\
84	84\\
85	85\\
86	86\\
87	87\\
88	88\\
89	89\\
90	90\\
91	91\\
92	92\\
93	93\\
94	94\\
95	95\\
96	96\\
97	97\\
98	98\\
99	99\\
100	100\\
101	101\\
102	102\\
103	103\\
104	104\\
105	105\\
106	106\\
107	107\\
108	108\\
109	109\\
110	110\\
111	111\\
112	112\\
113	113\\
114	114\\
115	115\\
116	116\\
117	117\\
118	118\\
119	119\\
120	120\\
121	121\\
122	122\\
123	123\\
124	124\\
125	125\\
126	126\\
127	127\\
128	128\\
129	129\\
130	130\\
131	131\\
132	132\\
133	133\\
134	134\\
135	135\\
136	136\\
137	137\\
138	138\\
139	139\\
140	140\\
141	141\\
142	142\\
143	143\\
144	144\\
145	145\\
146	146\\
147	147\\
148	148\\
149	149\\
150	150\\
151	151\\
152	152\\
153	153\\
154	154\\
155	155\\
156	156\\
157	157\\
158	158\\
159	159\\
160	160\\
161	161\\
162	162\\
163	163\\
164	164\\
165	165\\
166	166\\
167	167\\
168	168\\
169	169\\
170	170\\
};

\addplot [color=mycolor2, line width=0.5pt]
  table[row sep=crcr]{%
30	29.3185775747695\\
31	30.3425011917306\\
32	31.3661109424409\\
33	32.3894068269005\\
34	33.4123888451094\\
35	34.4350569970678\\
36	35.4574112827751\\
37	36.479451702232\\
38	37.5011782554382\\
39	38.5225909423936\\
40	39.5436897630985\\
41	40.5644747175524\\
42	41.5849458057559\\
43	42.6051030277086\\
44	43.6249463834105\\
45	44.6444758728617\\
46	45.6636914960622\\
47	46.6825932530122\\
48	47.7011811437115\\
49	48.7194551681597\\
50	49.7374153263574\\
51	50.7550616183047\\
52	51.772394044001\\
53	52.7894126034468\\
54	53.8061172966416\\
55	54.8225081235859\\
56	55.8385850842794\\
57	56.8543481787226\\
58	57.8697974069146\\
59	58.8849327688562\\
60	59.8997542645468\\
61	60.9142618939869\\
62	61.9284556571763\\
63	62.9423355541152\\
64	63.9559015848031\\
65	64.9691537492405\\
66	65.9820920474272\\
67	66.9947164793631\\
68	68.0070270450483\\
69	69.0190237444827\\
70	70.0307065776667\\
71	71.0420755445996\\
72	72.0531306452821\\
73	73.0638718797138\\
74	74.0742992478948\\
75	75.0844127498253\\
76	76.0942123855048\\
77	77.1036981549338\\
78	78.1128700581121\\
79	79.1217280950396\\
80	80.1302722657164\\
81	81.1385025701424\\
82	82.146419008318\\
83	83.1540215802429\\
84	84.1613102859169\\
85	85.1682851253403\\
86	86.1749460985129\\
87	87.1812932054347\\
88	88.1873264461061\\
89	89.1930458205267\\
90	90.1984513286966\\
91	91.2035429706158\\
92	92.2083207462842\\
93	93.2127846557018\\
94	94.216934698869\\
95	95.2207708757854\\
96	96.2242931864511\\
97	97.2275016308661\\
98	98.2303962090303\\
99	99.232976920944\\
100	100.235243766607\\
101	101.237196746019\\
102	102.23883585918\\
103	103.240161106091\\
104	104.241172486752\\
105	105.241870001161\\
106	106.24225364932\\
107	107.242323431228\\
108	108.242079346885\\
109	109.241521396291\\
110	110.240649579448\\
111	111.239463896353\\
112	112.237964347007\\
113	113.236150931411\\
114	114.234023649564\\
115	115.231582501467\\
116	116.228827487118\\
117	117.225758606519\\
118	118.22237585967\\
119	119.218679246569\\
120	120.214668767218\\
121	121.210344421616\\
122	122.205706209764\\
123	123.200754131661\\
124	124.195488187307\\
125	125.189908376703\\
126	126.184014699847\\
127	127.177807156741\\
128	128.171285747385\\
129	129.164450471777\\
130	130.157301329919\\
131	131.149838321811\\
132	132.142061447451\\
133	133.133970706841\\
134	134.12556609998\\
135	135.116847626869\\
136	136.107815287507\\
137	137.098469081894\\
138	138.08880901003\\
139	139.078835071915\\
140	140.068547267551\\
141	141.057945596935\\
142	142.047030060069\\
143	143.035800656951\\
144	144.024257387584\\
145	145.012400251965\\
146	146.000229250096\\
147	146.987744381976\\
148	147.974945647606\\
149	148.961833046984\\
150	149.948406580112\\
151	150.934666246989\\
152	151.920612047616\\
153	152.906243981992\\
154	153.891562050117\\
155	154.876566251992\\
156	155.861256587616\\
157	156.845633056989\\
158	157.829695660111\\
159	158.813444396983\\
160	159.796879267604\\
161	160.780000271974\\
162	161.762807410094\\
163	162.745300681963\\
164	163.727480087581\\
165	164.709345626949\\
166	165.690897300066\\
167	166.672135106932\\
168	167.653059047547\\
169	168.633669121912\\
170	169.613965330026\\
};
\nextgroupplot[%
height=1.7cm,
ylabel style={align=center},
ylabel={Error $[\%]$},
ymin=-0.5,
ymax=2.5,
]
\addplot [color=black]
  table[row sep=crcr]{%
30	2.27140808410174\\
31	2.12096389764331\\
32	1.98090330487213\\
33	1.85028234272572\\
34	1.72826810261952\\
35	1.61412286552068\\
36	1.50719088118024\\
37	1.40688729126483\\
38	1.31268880147852\\
39	1.22412578873445\\
40	1.14077559225372\\
41	1.06225678645753\\
42	0.988224272009845\\
43	0.918365051840528\\
44	0.852394583157883\\
45	0.790053615862816\\
46	0.731105443343076\\
47	0.675333504229361\\
48	0.622539283934476\\
49	0.572540473143547\\
50	0.525169347285157\\
51	0.480271336657379\\
52	0.437703761536558\\
53	0.397334710477725\\
54	0.359042043256358\\
55	0.322712502571133\\
56	0.288240920929564\\
57	0.255529511013052\\
58	0.224487229457537\\
59	0.195029205328398\\
60	0.167076225755321\\
61	0.140554272152582\\
62	0.115394101328559\\
63	0.0915308664837815\\
64	0.068903773745177\\
65	0.0474557703992236\\
66	0.0271332614739802\\
67	0.00788585169687788\\
68	-0.0103338897768753\\
69	-0.0275706441778186\\
70	-0.043866539523825\\
71	-0.0592613304219689\\
72	-0.07379256289176\\
73	-0.087495725635344\\
74	-0.100404389046988\\
75	-0.112550333100406\\
76	-0.123963665137881\\
77	-0.134672928485467\\
78	-0.144705202707813\\
79	-0.154086196252695\\
80	-0.1628403321455\\
81	-0.17099082733634\\
82	-0.178559766241488\\
83	-0.185568168967295\\
84	-0.192036054663024\\
85	-0.197982500400325\\
86	-0.203425695945189\\
87	-0.208382994752527\\
88	-0.212870961484187\\
89	-0.216905416322173\\
90	-0.220501476329578\\
91	-0.22367359408325\\
92	-0.226435593787129\\
93	-0.228800705055695\\
94	-0.230781594541471\\
95	-0.232390395563614\\
96	-0.233638735886584\\
97	-0.234537763779456\\
98	-0.235098172479858\\
99	-0.235330223175769\\
100	-0.235243766606715\\
101	-0.234848263385101\\
102	-0.234152803118093\\
103	-0.233166122418648\\
104	-0.231896621876451\\
105	-0.230352382057856\\
106	-0.228541178603363\\
107	-0.226470496474409\\
108	-0.224147543411951\\
109	-0.221579262652716\\
110	-0.218772344952343\\
111	-0.215733239957326\\
112	-0.212468166970746\\
113	-0.208983125142557\\
114	-0.205283903126485\\
115	-0.201376088231793\\
116	-0.197265075101868\\
117	-0.192956073948138\\
118	-0.188454118364186\\
119	-0.18376407274739\\
120	-0.178890639348595\\
121	-0.173838364972288\\
122	-0.168611647347695\\
123	-0.163214741187689\\
124	-0.157651763957309\\
125	-0.151926701362026\\
126	-0.146043412577185\\
127	-0.14000563522931\\
128	-0.13381699014432\\
129	-0.127480985873799\\
130	-0.121001023014786\\
131	-0.114380398328631\\
132	-0.107622308675021\\
133	-0.100729854767742\\
134	-0.0937060447612831\\
135	-0.0865537976805113\\
136	-0.0792759466959748\\
137	-0.0718752422581003\\
138	-0.0643543550941342\\
139	-0.0567158790758833\\
140	-0.048962333964734\\
141	-0.0410961680389966\\
142	-0.0331197606117153\\
143	-0.0250354244416043\\
144	-0.0168454080441983\\
145	-0.00855189790689988\\
146	-0.000157020613658442\\
147	0.00833715511827044\\
148	0.0169286164827164\\
149	0.0256154047085548\\
150	0.034395613258539\\
151	0.0432673860997577\\
152	0.052228916042028\\
153	0.0612784431424148\\
154	0.0704142531705653\\
155	0.0796346761343109\\
156	0.0889380848618132\\
157	0.0983228936376783\\
158	0.107787556891709\\
159	0.117330567935282\\
160	0.126950457747537\\
161	0.136645793804834\\
162	0.146415178954471\\
163	0.156257250329525\\
164	0.166170678304147\\
165	0.176154165485565\\
166	0.186206445743588\\
167	0.196326283274394\\
168	0.206512471698161\\
169	0.21676383318817\\
170	0.22707921763167\\
};
\end{groupplot}
\end{tikzpicture}%
	\caption{Semi-definite fit of the \gls{acr:em} temperature together with the relative error. The normalized RMSE is 0.15\% w.r.t. the maximum temperature.}
	\label{fig:EM_temp_fit}
\end{figure}
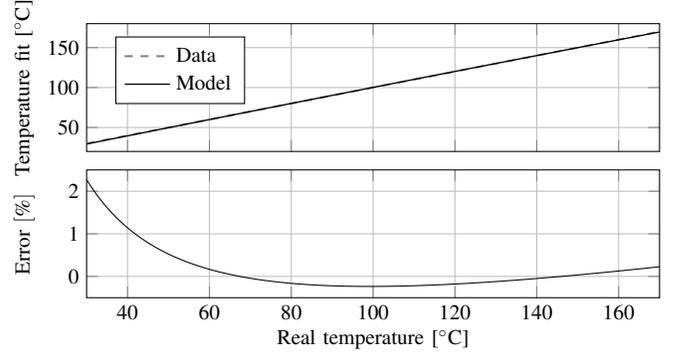

\nobreak\vspace{-5pt}
\begingroup
\allowdisplaybreaks
\begin{small}
	\begin{equation}
		\overline{\vartheta}_{\mathrm{m}}(s) = x_{\mathrm{m,\vartheta}}^{\top}(s) Q_{\mathrm{m,\vartheta}}x_{\mathrm{m,\vartheta}}(s) + \vartheta_\mathrm{m,0},
		\label{eq:EM_temp_fit}
	\end{equation}
\end{small}%
\endgroup
where $x_{\mathrm{m,\vartheta}}(s) = \left[1 \  \vartheta_\mathrm{m}(s) \right]^\top$, $\vartheta_\mathrm{m,0}$ is an offset required to obtain positive values and $Q_{\mathrm{m,\vartheta}} \in \mathbb{S}_-^2$ is a negative semi-definite matrix of coefficients, obtained through semi-definite programming. The semi-definite fit of the temperature together with the relative error is shown in Fig.~\ref{fig:EM_temp_fit}. We select a negative semi-definite matrix, since we require an upper bound on the \gls{acr:em} temperature. Translating~\eqref{eq:EM_cooling_time} to forces and substituting the offset of~\eqref{eq:EM_temp_fit} results in

\par\nobreak\vspace{-5pt}
\begingroup
\allowdisplaybreaks
\begin{small}
	\begin{equation}
		F_{\mathrm{m,c}}(s) = \revTCST{\left({ \overline{F}_{\mathrm{m,c}}(s) +\left(\vartheta_\mathrm{m,0} - \vartheta_\mathrm{m,c}\right)\cdot \frac{\mathrm{d}t}{\mathrm{d}s}(s)}\right)\cdot {\lambda_{\mathrm{m}}}},
		\label{eq:EM_cooling_lin}
	\end{equation}
\end{small}%
\endgroup 
where $\overline{F}_{\mathrm{m,c}}(s)$ is an intermediate variable used to obtain a convex formulation through
\par\nobreak\vspace{-5pt}
\begingroup
\allowdisplaybreaks
\begin{small}
	\begin{equation}
		\overline{F}_{\mathrm{m,c}}(s)\cdot v(s) \leq x_{\mathrm{m,\vartheta}}^{\top}(s) Q_{\mathrm{m,\vartheta}}x_{\mathrm{m,\vartheta}}(s),
		\label{eq:EM_cooling_cone}
	\end{equation}
\end{small}%
\endgroup
which can be written as a second-order conic constraint (see Appendix).
Since the cooling circuit consists of liquid coolant flowing through the \gls{acr:em} to a set of radiators, the coolant temperature cannot exceed the \gls{acr:em} temperature, resulting in 
\par\nobreak\vspace{-5pt}
\begingroup
\allowdisplaybreaks
\begin{small}
	\begin{equation}
		F_{\mathrm{m,c}}(s) \geq 0.
		\label{eq:EM_cooling_positive}
	\end{equation}
\end{small}%
\endgroup

\begin{figure}[!tb]
	\centering
	\input{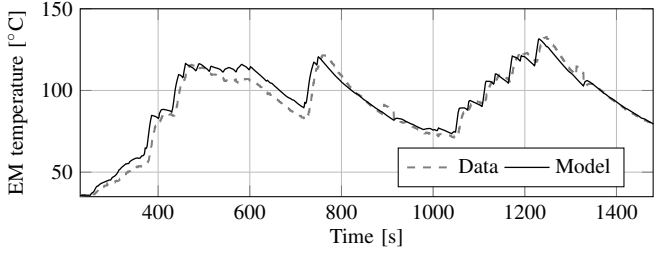}
	\caption{Comparison between the thermal model and one of the vehicle telemetry data samples. In total, three data sets were used to fit the model parameters, of which the average normalized RMSE was 4.12\% w.r.t. the maximum temperature of each data set.}
	\label{fig:EM_temp_fit_data}
\end{figure}

To obtain the thermal parameters of the \gls{acr:em}, we apply the convex loss- and cooling model to a combination of data sets recorded from vehicle telemetry. Fig.~\ref{fig:EM_temp_fit_data} shows a comparison between the \gls{acr:em} temperature from one of the data sets and the thermal model.

To prevent the \gls{acr:em} from overheating, we define an upper bound on the temperature through
\par\nobreak\vspace{-5pt}
\begingroup
\allowdisplaybreaks
\begin{small}
	\begin{equation}
		\vartheta_\mathrm{m}(s) \leq \vartheta_\mathrm{m,max},
		\label{eq:EM_temp_upper}
	\end{equation}
\end{small}%
\endgroup
where $\vartheta_\mathrm{m,max}$ is the maximum temperature of the \gls{acr:em}. 
Lastly, we specify an initial value for the \gls{acr:em} temperature as
\par\nobreak\vspace{-5pt}
\begingroup
\allowdisplaybreaks
\begin{small}
	\begin{equation}
		\vartheta_\mathrm{m}(0) = \vartheta_\mathrm{m,init},
		\label{eq:EM_temp_initial}
	\end{equation}
\end{small}%
\endgroup
where $\vartheta_\mathrm{m,init}$ is the initial value for the temperature and is calculated during pre-processing using a lookup table that has the charge time as an input. 

\subsection{Inverter}
In this section, we derive a quadratic loss model for the inverter losses. As opposed to the \gls{acr:em}, we do not model the inverter temperature, since we assume that the motor-inverter combination is designed such that the \gls{acr:em} is thermally limiting. We apply the general quadratic power loss model of the form
\par\nobreak\vspace{-5pt}
\begingroup
\allowdisplaybreaks
\begin{small}
	\begin{equation}
		P_{\mathrm{dc}}(s) = \alpha_\mathrm{inv} \cdot P_{\mathrm{ac}}^2(s)+P_{\mathrm{ac}}(s),
	\end{equation}
\end{small}%
\endgroup
where \revTCST{$\alpha_\mathrm{inv}\geq 0$} is an efficiency parameter, subject to identification. Converting this constraint to forces, rewriting and relaxing results in
\par\nobreak\vspace{-5pt}
\begingroup
\allowdisplaybreaks
\begin{small}
	\begin{equation} \label{eq:eff_inverter}
		(F_{\mathrm{dc}}(s)-F_{\mathrm{ac}}(s))\cdot \frac{\mathrm{d}t}{\mathrm{d}s}(s) \geq \alpha_\mathrm{inv} \cdot F_{\mathrm{ac}}^2(s),
	\end{equation}
\end{small}%
\endgroup
where $F_{\mathrm{dc}}(s)$ is the force equivalent to the electrical inverter power. The convexity of this constraint is shown in the Appendix by writing it as a second-order conic constraint. 
\par
\subsection{Battery}\label{sec:batt}
This section derives a model for the battery efficiency and the power-split between the electrical inverter power and the auxiliary component power. The latter can be observed from Fig.~\ref{fig:topology} and is written as
\par\nobreak\vspace{-5pt}
\begingroup
\allowdisplaybreaks
\begin{small}
	\begin{equation}
		P_{\mathrm{b}}(s) = P_{\mathrm{dc}}(s) + P_{\mathrm{aux}},
		\label{eq:P_bat_split}
	\end{equation}
\end{small}%
\endgroup  
where $P_{\mathrm{b}}(s)$ is the battery power at the terminals. Here, the auxiliary component supply is assumed to be constant and uni-directional, while the other powers are bi-directional. Converting \eqref{eq:P_bat_split} to forces results in
\par\nobreak\vspace{-5pt}
\begingroup
\allowdisplaybreaks
\begin{small}
	\begin{equation} \label{eq:Paux}
		F_{\mathrm{b}}(s) = F_{\mathrm{dc}}(s) +P_{\mathrm{aux}}\cdot \frac{\mathrm{d}t}{\mathrm{d}s}(s),
	\end{equation}
\end{small}%
\endgroup 
where $F_{\mathrm{b}}(s)$ is the force equivalent of the battery power at the terminals.

The battery efficiency is mostly determined by its internal resistance $R_{\mathrm{0}}$ and open-circuit voltage $V_{\mathrm{oc}}$. We derive the battery losses $P_{\mathrm{b,l}}(E_\mathrm{b},\vartheta_\mathrm{b},P_\mathrm{i})$ from a \revTCST{quasistatic zero-order Thévenin equivalent circuit model~\cite[Chapter 4.5]{GuzzellaSciarretta2007}, whereby we only consider the open-circuit voltage and internal resistance}. The battery losses are given by
\par\nobreak\vspace{-5pt}
\begingroup
\allowdisplaybreaks
\begin{small}
	\begin{equation}
		P_{\mathrm{b,l}}(E_\mathrm{b},\vartheta_\mathrm{b},P_\mathrm{i}) = \frac{1}{P_{\mathrm{sc}}(E_\mathrm{b},\vartheta_\mathrm{b})} \cdot P_{\mathrm{i}}^2(s),
		\label{eq:Pbloss}
	\end{equation}
\end{small}%
\endgroup 
where $P_{\mathrm{sc}}(E_\mathrm{b},\vartheta_\mathrm{b}) = \frac{V_{\mathrm{oc}}^2(E_\mathrm{b})}{R_{\mathrm{0}}(E_\mathrm{b},\vartheta_\mathrm{b})}$ is the short-circuit power~\cite{BorsboomFahdzyanaEtAl2020}. In reality, both the internal resistance and open-circuit voltage are a function of the battery temperature and energy. \revTCST{Since the battery reaches its thermal limit at the end of the charging process and has to be cooled during driving, it is operated at a relatively large temperature window. Therefore, the heat generation of the battery due to losses can vary significantly during driving.} However, since the influence of the battery temperature on the open-circuit voltage is rather small (less than $0.25\%$ \revTCST{according to manufacturer datasheets}) within the operating window in racing scenarios, we neglect the dependency of the open-circuit voltage on temperature~\cite{ZhangXiaEtAl2018,FarmannSauer2017}. To capture the impact of temperature on the internal resistance, we apply a correction factor inversely proportional to the battery temperature temperature~\cite{RosewaterCoppEtAl2019,Lebkowski2017} to obtain 
\par\nobreak\vspace{-5pt}
\begingroup
\allowdisplaybreaks
\begin{small}
	\begin{equation}
		P_{\mathrm{sc}}(E_\mathrm{b},\vartheta_\mathrm{b}) = 
		\frac{V_{\mathrm{oc}}^2(E_\mathrm{b})}{R_{\mathrm{0}}(E_\mathrm{b})}\cdot \frac{\vartheta_\mathrm{b}(s)}{\vartheta_\mathrm{ref}},
		\label{eq:Psc}
	\end{equation}
\end{small}%
\endgroup
where $\vartheta_\mathrm{ref}$ represents the reference temperature at which the battery data are measured. Similarly to the thermal \gls{acr:em} model, we fit the short-circuit power in a convex manner through
\par
\begin{figure}[!tb]
	\centering
	\input{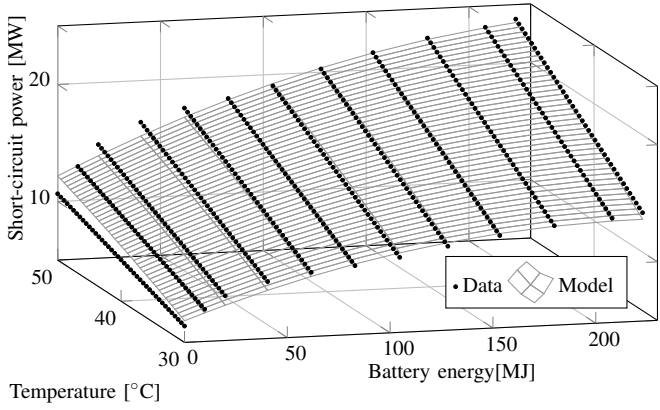}
	\caption{A temperature- and energy-dependent battery model \revTCST{compared to the calculated short-circuit power using manufacturer datasheets}. The normalized RMSE of the model is 2.23\% w.r.t. the maximum short-circuit power.}
	\label{fig:Bat_Psc_fit}
\end{figure}

\nobreak\vspace{-5pt}
\begingroup
\allowdisplaybreaks
\begin{small}
	\begin{equation}
		P_{\mathrm{sc}}(s) = x_{\mathrm{b,l}}^{\top}(s) Q_{\mathrm{b,l}}x_{\mathrm{b,l}}(s) + P_\mathrm{sc,0},
		\label{eq:Bat_Psc_fit}
	\end{equation}
\end{small}%
\endgroup
where $x_{\mathrm{b,l}}(s) = \left[1 \  E_\mathrm{b}(s) \  \vartheta_\mathrm{b}(s) \right]^\top$, $P_\mathrm{sc,0}$ is an offset required to obtain positive values and $Q_{\mathrm{b,l}} \in \mathbb{S}_-^2$ is a negative semi-definite matrix of coefficients, identified through semi-definite programming. Again, we select a negative semi-definite matrix, since it is optimal to maximize the short-circuit force, and thereby require an upper bound. The temperature- and energy-dependent model of the short-circuit power is shown in Fig.~\ref{fig:Bat_Psc_fit}. Translating~\eqref{eq:Psc} to forces and substituting the offset of~\eqref{eq:Bat_Psc_fit} results in

\par\nobreak\vspace{-5pt}
\begingroup
\allowdisplaybreaks
\begin{small}
	\begin{equation}
		F_{\mathrm{sc}}(s) =
		\overline{F}_\mathrm{sc}(s) + P_\mathrm{sc,0}\cdot \frac{\mathrm{d}t}{\mathrm{d}s}(s),
		\label{eq:Bat_Fsc_lin}
	\end{equation}
\end{small}%
\endgroup
where $F_{\mathrm{sc}}(s)$ is the short-circuit force and $\overline{F}_{\mathrm{sc}}(s)$ is an intermediate variable used to obtain a convex formulation through
\par\nobreak\vspace{-5pt}
\begingroup
\allowdisplaybreaks
\begin{small}
	\begin{equation}
		\overline{F}_{\mathrm{sc}}(s)\cdot v(s) \leq x_{\mathrm{b,l}}^{\top}(s) Q_{\mathrm{b,l}}x_{\mathrm{b,l}}(s).
		\label{eq:Bat_Fsc_cone}
	\end{equation}
\end{small}%
\endgroup

To obtain the battery losses during discharging $F_{\mathrm{b,l}}(s)$, we translate \eqref{eq:Pbloss} to forces and relax it, which results in
\par\nobreak\vspace{-5pt}
\begingroup
\allowdisplaybreaks
\begin{small}
	\begin{equation}\label{eq:eff_bat}
		F_{\mathrm{b,l}}(s)\cdot F_{\mathrm{sc}}(s) \geq F_{\mathrm{i}}^2(s),
	\end{equation}
\end{small}%
\endgroup
where $F_{\mathrm{i}}(s)$ is the internal battery force, which ultimately dictates a change in battery energy. To prevent the battery losses and short-circuit force from cooling the battery, we explicitly define 
\par\nobreak\vspace{-5pt}
\begingroup
\allowdisplaybreaks
\begin{small}
	\begin{align} 
		F_{\mathrm{b,l}}(s) &\geq 0, \label{eq:Bat_loss_positive}\\
		F_{\mathrm{sc}}(s) &\geq  0. \label{eq:Bat_Fsc_positive}
	\end{align} 
\end{small}%
\endgroup  

In general, the heat generation of a Li-ion battery is different between charging and discharging conditions. \revTCST{As the direction of power flow frequently changes during a lap, it is important to capture the additional heat generation during negative power flow.} To capture this behavior, we add an additional term to the battery losses to represent distinct heat generation during charging, compared to discharging. \revTCST{This means that we effectively have a different battery efficiency during charging and discharging~\cite{GuzzellaSciarretta2007}}. Since this term should only be present during negative power flow, we implement a set of inequality constraints similar to the final drive losses as
\par\nobreak\vspace{-5pt}
\begingroup
\allowdisplaybreaks
\begin{small}
	\begin{alignat}{2} 
		&F_{\mathrm{i}}(s) \geq \quad &&F_{\mathrm{b}}(s) + F_{\mathrm{b,l}}(s) , \label{eq:Bat_loss_disch}\\
		&F_{\mathrm{i}}(s) \geq  (1-\alpha_\mathrm{ch})\cdot &&F_{\mathrm{b}}(s) + F_{\mathrm{b,l}}(s) , \label{eq:Bat_loss_ch}
	\end{alignat} 
\end{small}%
\endgroup
	%
where $F_{\mathrm{b}}(s)$ is the battery force at the terminals and $\alpha_\mathrm{ch}$ is a coefficient that represents the distinct charging losses. Note that $\alpha_\mathrm{ch} \geq 0$ results in higher losses during charging, compared to discharging, whereas $\alpha_\mathrm{ch} \leq 0$ results in higher losses during discharging, compared to charging. In energy-limited scenarios, \eqref{eq:Bat_loss_disch} will hold with equality during discharging, whereas~\eqref{eq:Bat_loss_ch} will hold with equality during charging.


In contrast to the \gls{acr:em} cooling, where the difference between the \gls{acr:em} temperature and the ambient temperature (in the order of~$100^\circ$C) is sufficient to apply radiators, the difference between the battery temperature and the ambient air is significantly lower (in the order of~$25^\circ$C). Therefore, it is common to apply a refrigerant circuit instead of radiators to cool the battery during fast-charging pit stops and driving, which allows the coolant temperature to drop below ambient level. 
Again, all losses are assumed to be converted to heat, thereby changing the battery temperature according to the first-order temperature \gls{acr:ODE} given by
\par\nobreak\vspace{-5pt}
\begingroup
\allowdisplaybreaks
\begin{small}
	\begin{equation}
		C_{\mathrm{b}}\cdot \frac{\mathrm{d}\vartheta_{\mathrm{b}}}{\mathrm{d}t}(s) = P_{\mathrm{b,l}}(s) - P_{\mathrm{b,c}}(s),
		\label{eq:Bat_thermal_time}
	\end{equation}
\end{small}%
\endgroup
where $C_{\mathrm{b}}$ is the total lumped thermal capacity of the battery, $\vartheta_{\mathrm{b}}(s)$ is the temperature of the battery and $P_{\mathrm{b,c}}(s) \geq 0$ represents the power outflow from the battery cells to the cooling liquid. Since we consider a battery cooling circuit where the coolant temperature can be actively controlled, the cooling power is free within the bounds defined as
\par\nobreak\vspace{-5pt}
\begingroup
\allowdisplaybreaks
\begin{small}
	\begin{equation}
		0 \leq P_{\mathrm{b,c}}(s) \leq \revTCST{\left({\vartheta_{\mathrm{b}}(s)-\vartheta_{\mathrm{b,c}}}\right)\cdot {\lambda_{\mathrm{b}}}},
		\label{eq:Bat_cooling_time}
	\end{equation}
\end{small}%
\endgroup
where $\vartheta_{\mathrm{b,c}}$ represents the lowest achievable temperature of the cooling liquid and \revTCST{$\lambda_{\mathrm{b}}$ is the thermal conductivity} between the battery cells and the cooling liquid, where again we assume both parameters to be constant.  Rewriting~\eqref{eq:Bat_thermal_time} to space domain results in
\par\nobreak\vspace{-5pt}
\begingroup
\allowdisplaybreaks
\begin{small}
	\begin{equation}
		C_{\mathrm{b}}\cdot \frac{\mathrm{d}\vartheta_{\mathrm{b}}}{\mathrm{d}s}(s) = F_{\mathrm{i}}(s) - F_{\mathrm{b}}(s) - F_{\mathrm{b,c}}(s),
		\label{eq:Bat_thermal_space}
	\end{equation}
\end{small}%
\endgroup
where $ F_{\mathrm{b,c}}(s) \geq 0$ is the force-equivalence of the battery cooling power. Note that we explicitly use the difference between the internal battery force and the battery force at the terminals to include the additional charging losses. Similarly as with the \gls{acr:em} cooling, we approximate the battery temperature as
\par\nobreak\vspace{-5pt}
\begingroup
\allowdisplaybreaks
\begin{small}
	\begin{equation}
		\overline{\vartheta}_{\mathrm{b}}(s) = x_{\mathrm{b,\vartheta}}^{\top}(s) Q_{\mathrm{b,\vartheta}}x_{\mathrm{b,\vartheta}}(s) + \vartheta_\mathrm{b,0},
		\label{eq:Bat_temp_fit}
	\end{equation}
\end{small}%
\endgroup
where $x_{\mathrm{b,\vartheta}}(s) = \left[1 \  \vartheta_\mathrm{b}(s) \right]^\top$, $\vartheta_\mathrm{b,0}$ is an offset required to obtain positive values and $Q_{\mathrm{b,\vartheta}} \in \mathbb{S}_-^2$ is a negative semi-definite matrix of coefficients, obtained through semi-definite programming. The semi-definite fit of the battery temperature is similar to Fig.~\ref{fig:EM_temp_fit}, except that the normalized RMSE is reduced to $0.034\%$ due to the smaller temperature window. Translating~\eqref{eq:Bat_cooling_time} to forces and substituting the offset of~\eqref{eq:Bat_temp_fit} results in
\par\nobreak\vspace{-5pt}
\begingroup
\allowdisplaybreaks
\begin{small}
	\begin{equation}
		F_{\mathrm{b,c}}(s) = \revTCST{\left({ \overline{F}_{\mathrm{b,c}}(s) +\left(\vartheta_\mathrm{b,0} - \vartheta_\mathrm{b,c}\right)\cdot \frac{\mathrm{d}t}{\mathrm{d}s}(s)}\right)\cdot {\lambda_{\mathrm{b}}}},
		\label{eq:Bat_cooling_lin}
	\end{equation}
\end{small}%
\endgroup 
where $\overline{F}_{\mathrm{b,c}}(s)$ is an intermediate variable used to obtain a convex formulation (see Appendix) through
\par\nobreak\vspace{-5pt}
\begingroup
\allowdisplaybreaks
\begin{small}
	\begin{equation}
		\overline{F}_{\mathrm{b,c}}(s)\cdot v(s) \leq x_{\mathrm{b,\vartheta}}^{\top}(s) Q_{\mathrm{b,\vartheta}}x_{\mathrm{b,\vartheta}}(s).
		\label{eq:Bat_cooling_cone}
	\end{equation}
\end{small}%
\endgroup
As the coolant temperature cannot exceed the battery temperature, we define
\par\nobreak\vspace{-5pt}
\begingroup
\allowdisplaybreaks
\begin{small}
	\begin{equation}
		F_{\mathrm{b,c}}(s) \geq 0.
		\label{eq:Bat_cooling_positive}
	\end{equation}
\end{small}%
\endgroup

\begin{figure}[!t]
		\centering
		\input{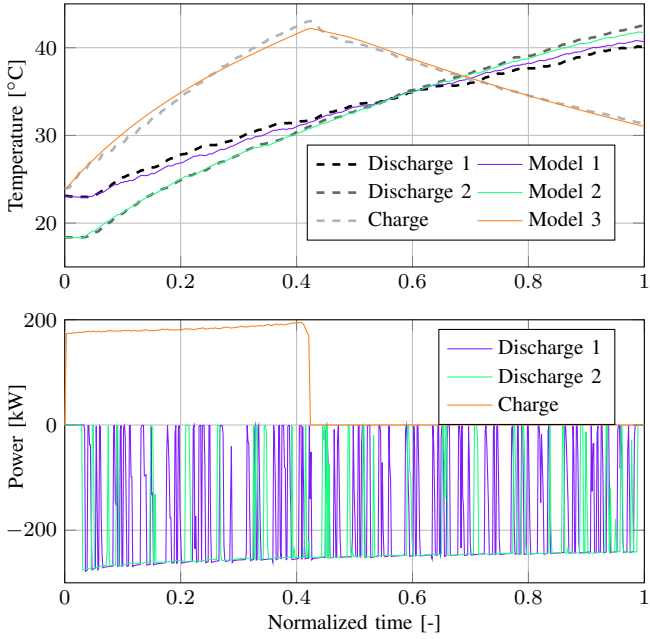}
		\caption{Comparison between the thermal model and three data sets, of which two represent a dynamic discharge (violet and green) and one represents a constant current charge profile followed by a cooldown period (orange). The average normalized RMSE is 0.87\% w.r.t. the maximum temperature of each data set. Both the modeled and the measured temperatures were obtained with the corresponding battery power profiles shown in the lower plot.}
		\label{fig:Bat_temp_fit_data}
\end{figure}

To obtain the thermal parameters of the battery, we apply the convex loss- and cooling model to a combination of data sets recorded from vehicle telemetry. Fig.~\ref{fig:Bat_temp_fit_data} shows a comparison between the battery temperature from various data sets and the thermal model, demonstrating that our model can accurately capture dynamic discharging as well as charging profiles. 

To ensure safe operation of the battery, we define an upper bound on the temperature through
\par\nobreak\vspace{-5pt}
\begingroup
\allowdisplaybreaks
\begin{small}
	\begin{equation}
		\vartheta_\mathrm{b}(s) \leq \vartheta_\mathrm{b,max},
		\label{eq:Bat_temp_upper}
	\end{equation}
\end{small}%
\endgroup
where $\vartheta_\mathrm{b,max}$ is the maximum temperature of the battery. Due to the relatively long charge time, compared to the refueling time of a conventional race car, it is essential to minimize the charge time. Therefore, we assume that the battery temperature reaches the upper limit at the end of charging, since the temperature is the main limitation. Hence, we enforce the initial battery temperature to be at the upper bound through
\par\nobreak\vspace{-5pt}
\begingroup
\allowdisplaybreaks
\begin{small}
	\begin{equation}
		\vartheta_\mathrm{b}(0) = \vartheta_\mathrm{b,max}.
		\label{eq:Bat_temp_initial}
	\end{equation}
\end{small}%
\endgroup
Lastly, we specify a terminal value for the battery temperature as
\par\nobreak\vspace{-5pt}
\begingroup
\allowdisplaybreaks
\begin{small}
	\begin{equation}
		\vartheta_\mathrm{b}(S_\mathrm{stint}) \leq \vartheta_\mathrm{b,N},
		\label{eq:EM_temp_final}
	\end{equation}
\end{small}%
\endgroup
where $\vartheta_\mathrm{b,N}$ is the terminal value for the temperature and is calculated during pre-processing using a lookup table having the charge time as an input.

The energy consumption of the battery is modeled as
\par\nobreak\vspace{-5pt}
\begingroup
\allowdisplaybreaks
\begin{small}
	\begin{equation} \label{eq:dEbat_ds}
		\frac{\mathrm{d}}{\mathrm{d}s}E_{\mathrm{b}}(s) = -F_{\mathrm{i}}(s),
	\end{equation}
\end{small}%
\endgroup
and we constrain the battery energy as
\par\nobreak\vspace{-5pt}
\begingroup
\allowdisplaybreaks
\begin{small}
	\begin{equation}
		E_{\mathrm{b,min}} \leq E_{\mathrm{b}}(s) \leq E_{\mathrm{b,max}},
		\label{eq:Ebat_bound}
	\end{equation} 
\end{small}%
\endgroup
where $E_{\mathrm{b,min}}$ and $E_{\mathrm{b,max}}$ correspond to the battery energy at the lower and upper \gls{acr:soe} bound, respectively. Since the voltage of the battery increases with the battery energy, it is optimal to operate the battery at higher energy levels. Therefore, we set the initial battery energy to the upper bound and constrain the terminal battery energy according to the amount of energy charged after driving as
\par\nobreak\vspace{-5pt}
\begingroup
\allowdisplaybreaks
\begin{small}
	\begin{align}
		E_{\mathrm{b}}(0) &= E_\mathrm{b,0}, \label{eq:Ebat_initial} \\  	
		E_{\mathrm{b}}(S_\mathrm{stint}) &\geq E_{\mathrm{b,0}} - E_{\mathrm{b,ch}},
		\label{eq:Ebat_final}
	\end{align} 
\end{small}%
\endgroup
where $E_\mathrm{b,0}$ is the initial battery energy and $E_{\mathrm{b,ch}}$ is the energy the battery receives during charging. In this way, the battery is guaranteed to be charged back to the upper energy bound in energy-limited scenarios. To calculate the battery energy during charging, we leverage a lookup table with input charge time $t_\mathrm{charge}$ and output $E_\mathrm{b,ch}$ for a given charging current profile during pre-processing. 

\subsection{Low-level Optimization Problem}
\label{sec:Methodology - low-level optimization}
This section presents the minimum-stint-time control problem of the electric race car. Given a predefined stint length and charge time we formulate the control problem using the state variables $x = (E_{\mathrm{kin}},E_{\mathrm{b}}, \vartheta_\mathrm{b}, \vartheta_\mathrm{m})$ and the control variables $u = (F_{\mathrm{m}},F_{\mathrm{brake,F}}, F_{\mathrm{brake,R}} )$ as follows:\\

\begin{prob}[Minimum-stint-time Control Strategy]\label{prob:low-level}
	The minimum-stint-time control strategies are the solution of
	\par\nobreak\vspace{-5pt}
	\begingroup
	\allowdisplaybreaks
	\begin{small}
		\begin{equation*}
			\begin{aligned}
				&\min \int_{0}^{S_{\mathrm{stint}}} {\dtds(s)}\,\mathrm{d}s ,\\
				&\textnormal{s.t. }  \eqref{eq:lethargy} -\eqref{eq:P_em,max}, \eqref{eq:eff_EM}, \eqref{eq:EM_thermal_space}, \eqref{eq:EM_cooling_lin} - \eqref{eq:EM_temp_initial}, \eqref{eq:eff_inverter}, \\
				& \qquad \eqref{eq:Paux}, \eqref{eq:Bat_Fsc_lin}-\eqref{eq:Bat_loss_ch}, \eqref{eq:Bat_thermal_space}, \eqref{eq:Bat_cooling_lin}-\eqref{eq:Ebat_final}.\\
			\end{aligned}
		\end{equation*}
	\end{small}%
	\endgroup
\end{prob}
\noindent Since the feasible domain and the cost function are convex, the low-level control problem is fully convex and therefore we can compute a globally optimal solution with standard nonlinear programming methods.

\section{High-level Race Optimization} \label{sec:high level}
In this section, we present the high-level maximum-race-distance control problem. First, we formulate the maximum-race-distance control problem that optimizes the stint length and charge time for a pre-defined number of pit stops. Second, we model the minimum stint time by leveraging the low-level control problem and optimizing for various combinations of stint length and initial battery energy, as was shown in Fig.~\ref{fig:framework}.
Finally, we extend the maximum-race-distance control problem to allow joint optimization of the stint length, charge time, and number of pit stops. 

\subsection{Mixed-integer Control Problem}
We define the high-level control problem for a pre-defined number of pit stops in \emph{stint domain}, so that we have a fixed and finite optimization horizon. Here, each index in the optimization variables represents a stint. The goal is then to maximize the driven distance as the sum of all completed laps during the stints as
\par\nobreak\vspace{-5pt}
\begingroup
\allowdisplaybreaks
\begin{small}
	\begin{equation}
		\max S_{\mathrm{race}} = \max \sum_{k=0}^{n_{\mathrm{stops}}} S_{\mathrm{lap}}\cdot N_{\mathrm{laps}}(k),
	\end{equation}
\end{small}%
\endgroup
where $S_{\mathrm{race}}$ is the total race distance, $n_{\mathrm{stops}}$ is the pre-defined number of pit stops, $N_{\mathrm{laps}}(k)\in\sN, \  \forall~k \in \left[0,...,~n_{\mathrm{stops}}~-~1\right]$ is the stint length and $\sN$ the set of natural numbers, and $S_{\mathrm{lap}}$ is the length of one lap. Since the vehicle starts and stops at the pit box, the stint length should be an integer number of laps. As it is unlikely that the vehicle is exactly at the finish line when the race time limit is reached, we allow the final stint length to be a non-integer number of laps. This way, we have $n_{\mathrm{stops}}+1$ stints for $n_\mathrm{stops}$ pit stops and thus we have $n_\mathrm{stops}$ integer stint lengths and one final non-integer stint length.

The race can be divided into the car driving a stint and the car recharging the battery during pit-stops. The total elapsed time $t_\mathrm{tot}(k)$ is then updated after every pit stop by adding both the time to complete the stint $t_{\mathrm{stint}}(k) \geq 0$ and the time spent on charging $t_{\mathrm{charge}}(k)\geq 0$ according to
\par\nobreak\vspace{-5pt}
\begingroup
\allowdisplaybreaks
\begin{small}
	\begin{equation} \label{eq:t_elapsed}
		t_{\mathrm{tot}}(k+1) = t_\mathrm{tot}(k) + t_{\mathrm{stint}}(k) + t_{\mathrm{charge}}(k), \quad \forall k \in [0,n_\mathrm{stops}-1].
	\end{equation}
\end{small}%
\endgroup
Since we do not have a pit stop after the final stint, the elapsed time after the final stint is defined as
\par\nobreak\vspace{-5pt}
\begingroup
\allowdisplaybreaks
\begin{small}
	\begin{equation} \label{eq:t_elapsed_final}
		t_{\mathrm{tot}}(n_\mathrm{stops}+1) = t_\mathrm{tot}(n_\mathrm{stops}) + t_{\mathrm{stint}}(n_\mathrm{stops}).
	\end{equation}
\end{small}%
\endgroup
Moreover, we define the elapsed time counter to start at the beginning of the race and the elapsed time after the final stint should not exceed the total race time limit $t_\mathrm{race}$, resulting in
\par\nobreak\vspace{-5pt}
\begingroup
\allowdisplaybreaks
\begin{small}
	\begin{align}
		t_\mathrm{tot}(0)&=0, \label{eq:tel_initial} \\  	
		t_\mathrm{tot}(n_\mathrm{stops}+1)&\leq t_\mathrm{race},
		\label{eq:t_race}
	\end{align} 
\end{small}%
\endgroup

We then decompose the total race into blocks consisting of the vehicle first driving a stint followed by a pit stop in which the battery is charged. Assuming that a stint is always energy-limited, the charge time uniquely defines the terminal battery energy for the prior stint and is not influenced by other stints. 
Furthermore, we assume that the battery is thermally-limited during charging, which allows us to pre-calculate the maximum terminal battery temperature for the prior stint through backwards integration of the battery temperature dynamics during charging.
Thereby, we uniquely base the terminal battery temperature on the charge time, without being influenced by other stints. 
Lastly, we assume that the \gls{acr:em} temperature reaches the upper limit at the end of driving and assume that the charge times across consecutive stints remains constant. This allows us to calculate the \gls{acr:em} temperature at the beginning of the stint by integrating the \gls{acr:em} temperature dynamics during charging. This way, we uniquely define the initial \gls{acr:em} temperature on the charge time. 

To ensure that the battery is not overcharged, we apply an upper bound on the charge time through
\par\nobreak\vspace{-5pt}
\begingroup
\allowdisplaybreaks
\begin{small}
	\begin{equation}\label{eq:t_charge_max}
		t_\mathrm{charge}(k) \leq t_{\mathrm{charge,max}},
	\end{equation}
\end{small}%
\endgroup
where $t_{\mathrm{charge,max}}$ is the maximum charge time corresponding to charging the battery from the lower to the upper energy level.
Finally, the time to complete the stint is obtained by solving the low-level control problem, which we explain in the subsequent section.

\subsection{Stint Time Model}\label{sec:stinttime}
\begin{figure}
	\centering
	\input{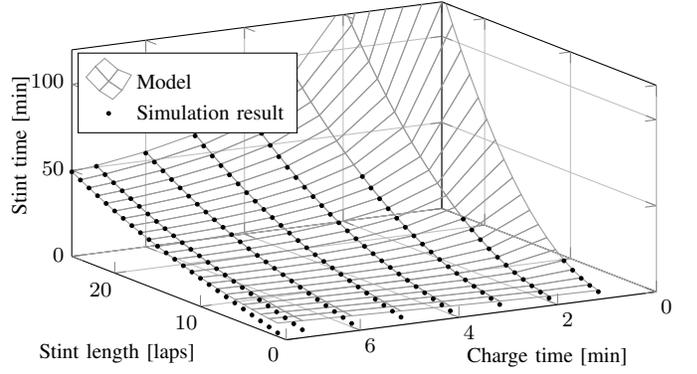}
	\caption{Fit of optimization results for a combination of stint lengths and charge times. The normalized RMSE of the fit is 2.2\% w.r.t. the maximum stint time.}
	\label{fig:t_stint_lut}
\end{figure}

In this section, we derive a method for modeling the stint time as a function of the stint length and charge time during the pit stop prior to the stint.
We solve the low-level control Problem~\ref{prob:low-level} for a combination of stint lengths and charge times to obtain the respective achievable minimum stint time.
In such manner, we can create the lookup table with stint time as a function of stint length and charge time shown in Fig.~\ref{fig:t_stint_lut}. Thereby, the charge time and terminal battery energy are linked through a pre-defined charging profile, cf. Section~\ref{sec:batt}, whereas the terminal battery temperature and initial \gls{acr:em} temperature are calculated during pre-processing based on the charge time solely, cf. Section~\ref{sec:EM},~\ref{sec:batt}.
As the stint time increases for larger stint lengths and shorter charge times, similar to the \gls{acr:em} loss fit in Section~\ref{sec:EM} above, we approximate the low-level optimization results via the continuous function
\par\nobreak\vspace{-5pt}
\begingroup
\allowdisplaybreaks
\begin{small}
	\begin{equation}
		t_{\mathrm{stint}}(k) = x_{\mathrm{s}}^{\top}(k) Q_{\mathrm{s}}x_{\mathrm{s}}(k),
		\label{eq:t_stint}
	\end{equation}
\end{small}%
\endgroup
where $x_{\mathrm{s}}(k) = \left[\frac{1}{\sqrt{t_{\mathrm{charge}}(k)}} \  \sqrt{t_{\mathrm{charge}}(k)} \  \frac{N_{\mathrm{laps}}(k)}{\sqrt{t_{\mathrm{charge}}(k)}}\right]^\top$ and $Q_{\mathrm{s}} \in \mathbb{S}_+^3$ is a symmetric positive semi-definite matrix of coefficients. The result of the fit is shown in Fig.~\ref{fig:t_stint_lut}. For a convex implementation, we relax and rewrite \eqref{eq:t_stint} to
\par\nobreak\vspace{-5pt}
\begingroup
\allowdisplaybreaks
\begin{small}
	\begin{equation}
		t_{\mathrm{stint}}(k)\cdot t_{\mathrm{charge}}(k) \geq y_{\mathrm{s}}(k)^\top Q_{\mathrm{s}}y_{\mathrm{s}}(k),
		\label{eq:t_stint2}
	\end{equation}
\end{small}%
\endgroup
where $y_{\mathrm{s}}(k) = \left[1 \ t_{\mathrm{charge}}(k) \  N_{\mathrm{laps}}(k)\right]^\top$, and convert this relaxation to a conic constraint~\cite{BoydVandenberghe2004} as
\par\nobreak\vspace{-5pt}
\begingroup
\allowdisplaybreaks
\begin{small}
	\begin{equation}
		t_{\mathrm{stint}}(k) + t_{\mathrm{charge}}(k) \geq 
		\begin{Vmatrix}
			2\cdot z_{\mathrm{s}}(k) \\
			t_{\mathrm{stint}}(k) - t_{\mathrm{charge}}(k)\\	
		\end{Vmatrix}_2,
		\label{eq:t_stint_cone}
	\end{equation}
\end{small}%
\endgroup
where $z_{\mathrm{s}}=L_{\mathrm{s}}y_{\mathrm{s}}(k)$ with $L_{\mathrm{s}}$ being the Cholesky factorization of $Q_{\mathrm{s}}$~\cite{BoydVandenberghe2004}. Since it is optimal to minimize stint time, this constraint will hold with equality at the optimum.

The final stint of the race is not followed by a pit stop in which the battery is charged, which means that the battery can be fully depleted and there is no margin needed in the battery temperature. Therefore, we separately model the final stint by solving the low-level control Problem \ref{prob:low-level} for a range of stint lengths, with a fixed charge time $t_\mathrm{charge}(n_\mathrm{stops})=t_\mathrm{charge,max}$ and terminal battery temperature $\vartheta_\mathrm{b,N}=\vartheta_\mathrm{b,max}$. With the charge time being fixed, we can then model the final stint time by a quadratic function with the stint length as 
\par
\begin{figure}
	\centering
%
%
\definecolor{mycolor1}{rgb}{0.50000,0.500,0.500}%
\definecolor{mycolor2}{rgb}{0.000,0.000,1.00}%
\definecolor{mycolor3}{rgb}{0.400,0.400,0.400}%
\definecolor{mycolor4}{rgb}{0.00,1.00,0.00}%
\definecolor{mycolor5}{rgb}{0.700,0.700,0.700}%
\definecolor{mycolor6}{rgb}{1.00,0.00,0.00}%
\begin{tikzpicture}
	\footnotesize
	\begin{axis}[%
		width=0.87\columnwidth,
		height = 2.5cm,
		at={(0,0)},
		scale only axis,
		xmin=0,
		xmax=25,
		ymin=0,
		ymax=50,
		y label style={yshift=-0.6cm},
		axis background/.style={fill=white},
		xmajorgrids,
		ymajorgrids,
		xlabel={Stint length [laps]},
		ylabel={Stint time [min]},
		x label style={yshift=0.2cm},
		legend style={at={(0.95,0.1)}, anchor=south east, legend cell align=left, align=left},
		legend entries={Simulation result, Model},
		legend columns=2
		]		
\addplot [color=black, only marks,mark size = 0.7pt,  mark=*, mark options={solid, black}]
table[row sep=crcr]{%
	1	1.80636600007539\\
	2	3.3566684262339\\
	3	4.94360483471172\\
	4	6.52414991364911\\
	5	8.10864607337751\\
	6	9.69187759269445\\
	7	11.2934003618537\\
	8	12.9197205935645\\
	9	14.5842328970786\\
	10	16.2912513478259\\
	11	18.015986000658\\
	12	19.7727825012918\\
	13	21.5642344994843\\
	14	23.3714921987263\\
	15	25.2498026992757\\
	16	27.1617532678081\\
	17	29.2240149577822\\
	18	31.4313498355705\\
	19	33.7109628309718\\
	20	36.0724439821621\\
	21	38.5121517724957\\
	22	41.0078904035309\\
	23	43.6337833373452\\
	24	46.3101482804768\\
	25	49.042345473629\\
};
\addplot [color=mycolor1]
table[row sep=crcr]{%
	1	2.38278458837538\\
	2	3.67576367294123\\
	3	5.02358524446496\\
	4	6.42624930294657\\
	5	7.88375584838608\\
	6	9.39610488078346\\
	7	10.9632964001387\\
	8	12.5853304064519\\
	9	14.2622068997229\\
	10	15.9939258799518\\
	11	17.7804873471387\\
	12	19.6218913012833\\
	13	21.5181377423859\\
	14	23.4692266704464\\
	15	25.4751580854647\\
	16	27.535931987441\\
	17	29.6515483763751\\
	18	31.8220072522671\\
	19	34.047308615117\\
	20	36.3274524649247\\
	21	38.6624388016904\\
	22	41.0522676254139\\
	23	43.4969389360953\\
	24	45.9964527337346\\
	25	48.5508090183318\\
};
\end{axis}
\end{tikzpicture}%
	\caption{Fit of the final stint times for a range of stint lengths. The normalized RMSE of the fit is \revTCST{0.60\%} w.r.t. the maximum stint time.}
	\label{fig:t_stint_final}
\end{figure}
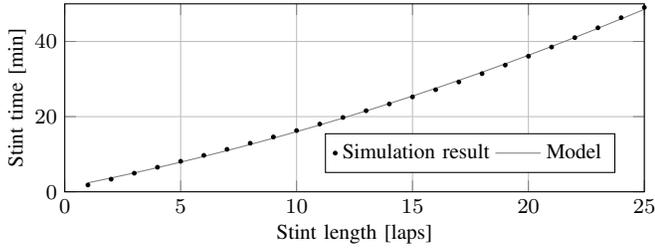

\nobreak\vspace{-5pt}
\begingroup
\allowdisplaybreaks
\begin{small}
	\begin{equation}
		t_{\mathrm{stint}}(n_\mathrm{stops}) \geq D_\mathrm{s,f}^\top x_\mathrm{s,f},
		\label{eq:t_stint_final}
	\end{equation}
\end{small}%
\endgroup
where $D_\mathrm{s,f}$ is a vector of coefficients and $x_\mathrm{s,f} = [N_\mathrm{laps}^2(n_\mathrm{stops}) \  N_\mathrm{laps}(n_\mathrm{stops}) \ 1]^\top$. Fig.~\ref{fig:t_stint_final} shows the quadratic fit of the final stint time model. 

\par
\subsection{Optimal Pit Stop Strategy}
In the previous sections, we introduced the objective and constraints for the high-level control problem when optimizing the race strategy for a pre-defined number of pit stops. In this section, we apply some modifications in order to jointly optimize the stint lengths, charge times and number of pit stops, thereby removing the need to search over a large space of pre-defined number of pit stops.

We define a binary variable $b_{\mathrm{pit}}(k)$ that indicates whether pit stop and stint $k$ is taken or skipped as
\par\nobreak\vspace{-5pt}
\begingroup
\allowdisplaybreaks
\begin{small}
	\begin{equation}
		b_{\mathrm{pit}}(k) =
		\begin{cases}
			0, & \text{if stint and stop skipped} \\
			1, & \text{if stint and stop taken}, \\
		\end{cases}
	\end{equation}
\end{small}%
\endgroup
and include it in \eqref{eq:t_stint} via the big-M formulation~\cite{RichardsHow2005}
\par\nobreak\vspace{-5pt}
\begingroup
\allowdisplaybreaks
\begin{small}
	\begin{equation}
		t_{\mathrm{stint}}(k) \geq x_{\mathrm{s}}(k)^\top Q_{\mathrm{s}}x_{\mathrm{s}}(k) - M\cdot (1-b_{\mathrm{pit}}(k)),
		\label{eq:t_stint_M}
	\end{equation}
\end{small}%
\endgroup
where $M \gg t_{\mathrm{stint,max}}$. This way, we obtain the original constraint if $b_{\mathrm{pit}}(k)=1$ and we obtain a negative lower bound when $b_{\mathrm{pit}}(k)=0$. By defining
\par\nobreak\vspace{-5pt}
\begingroup
\allowdisplaybreaks
\begin{small}
	\begin{align}
		t_{\mathrm{stint}}(k) &\geq 0, 	\label{eq:t_stint_pos}\\
		t_\mathrm{charge}(k) &\geq 0,
		\label{eq:t_charge_pos}
	\end{align}
\end{small}%
\endgroup
the $k$-th stint time and charge time will be pushed to zero, hence skipping the stint.
We convert \eqref{eq:t_stint_M} to a cone as
\par\nobreak\vspace{-5pt}
\begingroup
\allowdisplaybreaks
\begin{small}
	\begin{equation} \label{eq:t_stint_cone_M}
		\begin{split}
			M\cdot (1-b_{\mathrm{pit}}(k)) + t_{\mathrm{stint}}(k) + t_{\mathrm{charge}}(k)  \geq  \\
			\begin{Vmatrix}
				2\cdot z_{\mathrm{s}}(k) \\
				M\cdot (1-b_{\mathrm{pit}}(k)) +t_{\mathrm{stint}}(k) - t_{\mathrm{charge}}(k)   \\	
			\end{Vmatrix}_2.\\
		\end{split}
	\end{equation}
\end{small}%
\endgroup
Hence, whenever a stint is skipped, the corresponding stint time and charge time will be zero if an optimal solution is obtained. To prevent the stint length from diverging to infinity whenever the stint is actually skipped, i.e., $b_\mathrm{pit}(k) = 0$, we define an upper bound on stint length as
\par\nobreak\vspace{-5pt}
\begingroup
\allowdisplaybreaks
\begin{small}
	\begin{equation} \label{eq:Nlaps_M}
		N_{\mathrm{laps}}(k) \leq N_{\mathrm{laps,max}} \cdot b_{\mathrm{pit}} (k),
	\end{equation}
\end{small}%
\endgroup
where $N_\mathrm{laps,max}$ is the maximum stint length that was used to obtain the lookup table. This will ensure $N_{\mathrm{laps}}(k)=0$ whenever $b_{\mathrm{pit}}(k)=0$.
Since the final stint is not constrained by a lower terminal battery temperature and can deplete the battery without spending time on charging afterwards, the final stint has to be part of the optimal race strategy. We ensure this by writing 
\par\nobreak\vspace{-5pt}
\begingroup
\allowdisplaybreaks
\begin{small}
	\begin{equation} \label{eq:b_pit}
		b_{\mathrm{pit}}(k+1) \geq b_{\mathrm{pit}}(k),  \quad  \forall k \in [0,n_{\mathrm{stops}}-1].
	\end{equation}
\end{small}%
\endgroup
Suppose we set the length of the variables to $N$ and that $x$ stints and stops are skipped. Then the first $x$ entries in $b_\mathrm{pit}$ will be zero and the last $N-x$ entries will be one. In this way, the final stint is always taken and the search space for the solver is reduced.

\subsection{High-level Optimization Problem}
\label{sec:Methodology - high-level optimization}
This section presents the maximum-race-distance control problem of the electric race car. Given a predefined race time we formulate the control problem using the state variables $x=(t_\mathrm{tot})$ and the control variables $u=(t_{\mathrm{charge}},N_{\mathrm{laps}},b_{\mathrm{pit}})$ as follows:\\

\begin{prob}[Maximum-race-distance Strategies]\label{prob:high-level}
	The maximum-race-distance strategies are the solution of
	\par\nobreak\vspace{-5pt}
	\begingroup
	\allowdisplaybreaks
	\begin{small}
		\begin{equation*}
			\begin{aligned}
				& \max \sum_{k=0}^{n_{\mathrm{stops}}} S_{\mathrm{lap}}\cdot N_{\mathrm{laps}}(k) ,\\
				&\textnormal{s.t. }  \eqref{eq:t_elapsed}-\eqref{eq:t_charge_max}, \eqref{eq:t_stint_final}, \eqref{eq:t_stint_pos}-\eqref{eq:b_pit}. \\
			\end{aligned}
		\end{equation*}
	\end{small}%
	\endgroup
\end{prob}

\noindent Since Problem 2 can be solved with mixed-integer second-order conic programming solvers, we can guarantee global optimality upon convergence~\cite{Lee2012,BelottiKirchesEtAl2013}.

\section{Discussion}\label{sec:discussion}
A few comments are in order.
First, we assume that endurance racing tires do not degrade significantly within a single stint\revTCST{, considering that the vehicle is operated less at the grip limits and the stints are relatively short, compared to conventional endurance race cars. The tires can be changed every stint without additional time loss in the pit lane, since the charging time is much longer than the tire change time. Furthermore, we assume that the tires are operated at their average operating temperature.} Therefore, we disregard the influence of tire wear and tire temperature on the available grip level. Yet the high-level control problem can be readily extended to capture these dynamics if the lookup table is devised accounting for tire degradation. 
Second, we assume that the time gained from starting the race from the grid compared to the pit lane is negligible on a full endurance race. Hence, we do not separately optimize the first stint. 
Third, when the battery temperature is not an active constraint, the battery cooling relaxes in order to maximize the battery efficiency. Yet this can be interpreted as the battery coolant temperature being controlled by the refrigerant cooling system, which allows reduced cooling power, assuming that the system can cope with the requested coolant temperature changes.   
Fourth, in scenarios where the battery or \gls{acr:em} temperature is very limiting and the lower battery energy limit is not reached, it may occur that~\eqref{eq:lethargy} relaxes due to the positive contribution of a higher lethargy in~\eqref{eq:EM_cooling_lin},~\eqref{eq:Bat_Fsc_lin} and~\eqref{eq:Bat_cooling_lin}. In these cases, the solution will not reflect reality anymore and it should be disregarded. Yet, from a high-level perspective, these solutions cannot be part of the optimal race strategy, since charging the battery more than needed is sub-optimal. In fact, these situations were not found to be part of the optimal strategies in this study, as shown in Section~\ref{Results} below. 
If a stint with a relaxed lethargy is not part of the optimal race strategy, it is guaranteed that the physically correct solution of the same stint will not be part of the optimal strategy either, since it is always worse or equal, compared to the convex solution.
\revTCST{Fifth}, we assume a fixed value for the terminal \gls{acr:em} temperature when pre-calculating the bound on the initial \gls{acr:em} temperature. However, it can occur that this value is not reached, e.g., in scenarios where the \gls{acr:em} temperature is not an active constraint. In these situations, the temperature trajectories might not reflect reality, yet this does not affect the resulting solution, as there are no other states that depend on the \gls{acr:em} temperature.
\revTCST{Finally, we disregarded dynamic events, such as competitor interaction. Therefore, the results represent the theoretical optimal race strategy, which can be regarded as an upper bound on what is achievable in practice. The methods could be modified and applied online to account for such events, but this is beyond the current scope.}

\section{Results } \label{Results}
This section presents numerical results for both the low- and high-level control problem. We base our use case on the driven electric endurance race car of InMotion~\cite{InMotion}, shown in Fig.~\ref{fig:inmotion-car}, performing an 11 lap stint at the Zandvoort circuit for the low-level control problem and a \unit[6]{h} race at the same circuit for the high-level control problem. First, we discuss the numerical solutions for both control problems. Second, we validate the high-level control problem by comparing the optimal race strategy against fixed-pit-stop-number strategies and compare the results to the expected optimal combinations of stint length and charge time. \revTCST{Although the model is capable of capturing inclination and banking of the track, the results were generated for a flat track, since there is no open-access data of the Zandvoort circuit.} 
\par 
For the discretization of the model, we apply the trapezoidal method with a fixed step-size of $\Delta s = \unit[4]{m}$. We parse the low-level control problem with CasADi~\cite{AnderssonGillisEtAl2019} and solve it using IPOPT~\cite{WachterBiegler2006} combined with the MA57 linear solver~\cite{HSL}, whilst we parse the high-level control problem with YALMIP~\cite{Loefberg2004} and solve it using MOSEK~\cite{ApS2017}.
We perform the numerical optimization on an Intel Core i7-4710MQ \unit[2.5]{GHz} processor and \unit[8]{GB} of RAM. Thereby, the computation time for solving the low-level problem was about \unit[0.68]{s} of parsing and \unit[61]{s} of solving, whereas the high-level problem needed \unit[0.21]{s} of parsing and \unit[4.5]{s} of solving.

\begin{figure}[!tb]
	\centering
	\includegraphics[width=1\columnwidth,trim={8 0 17 0},clip]{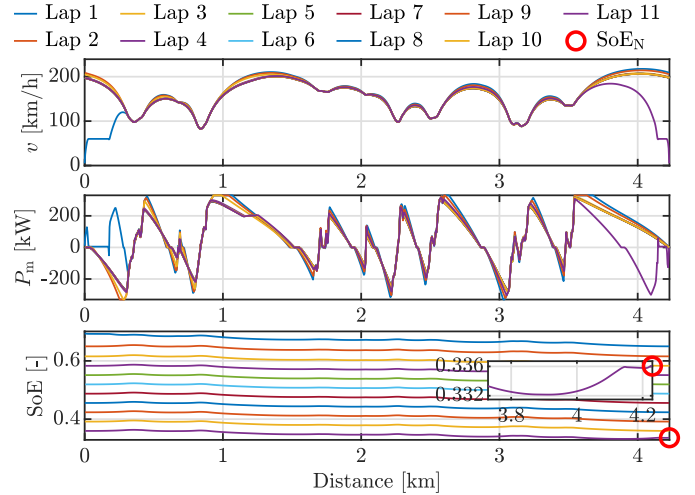}
	\caption{Velocity, EM power and battery SoE trajectories per lap for an 11 lap stint. The battery energy is an active constraint, thus the stint is energy-limited. The \gls{acr:em} power shows a gradual decrease at high velocities, thus indicating energy management. The zoom window corresponds to the final 500 m of the stint.}
	\label{fig:lap_results}
\end{figure}

\subsection{Low-level Optimization}\label{sec:results - low-level}
In this section, we compute the optimal trajectories for a stint of 11 laps around the Zandvoort circuit. We set the terminal battery capacity to the energy level corresponding to a \unit[5]{min} charge time using constant current charging, which means that the battery is partially charged. The total stint time is about \revTCST{ \unit[1126]{s}} with an average flying lap time of \unit[101]{s} (not counting the first and last lap).   
\par
\begin{figure}[!tb]
	\centering
	\input{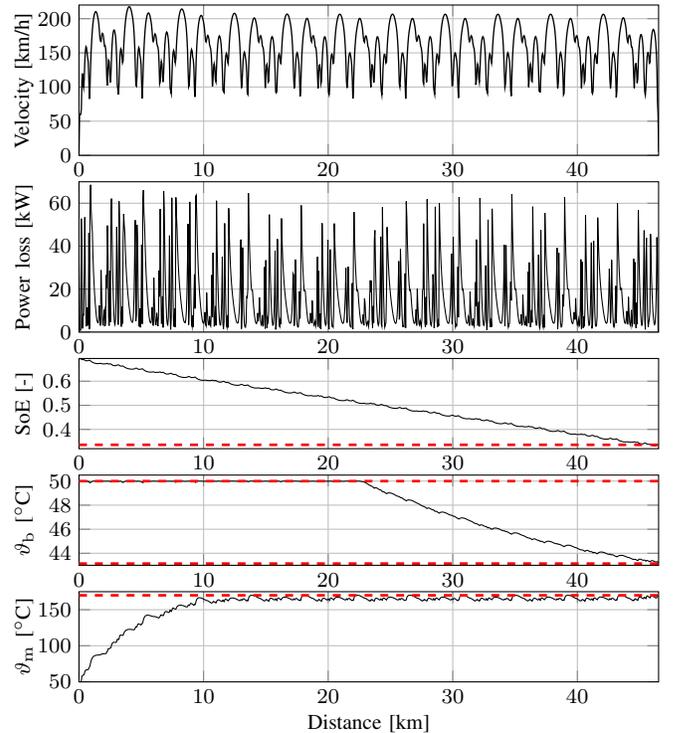}
	\caption{Velocity, powertrain losses, battery SoE, battery temperature and \gls{acr:em} temperature trajectories for an 11 lap stint. The battery energy is an active constraint, thus the stint is energy-limited and energy management is needed. Furthermore, the maximum \gls{acr:em} temperature is reached multiple times throughout the stint. Therefore, the velocity profile decreases when the maximum \gls{acr:em} temperature is reached.}
	\label{fig:Low_level_results}
\end{figure}

 
The velocity profile together with the EM power and \gls{acr:soe} per lap is shown in Fig.~\ref{fig:lap_results}. Furthermore, the total stint velocity profile together with the powertrain losses, \gls{acr:soe}, battery temperature and \gls{acr:em} temperature is shown in Fig.~\ref{fig:Low_level_results}.

First, we observe that the velocity profiles of the first four consecutive laps are slightly different, which is due to the \gls{acr:em} temperature being an active constraint. Since the \gls{acr:em} temperature starts below the limit, the velocity in the first laps is the highest and decreases as the \gls{acr:em} temperature reaches the upper limit. In this scenario, the terminal battery temperature is not an active constraint, since the battery temperature is kept at the maximum level during the first half of the stint by reducing the cooling power, thereby increasing the battery efficiency. In the second half of the stint, the battery temperature decreases gradually by applying maximum cooling to reach the required terminal value. If the terminal battery temperature were an active constraint, we would observe a gradual decrease in temperature over the entire stint, \revTCST{as shown in the Appendix for illustration purposes.}  
Second, the \gls{acr:em} power decreases gradually before the vehicle reaches a corner and regenerative braking is applied. However, both the power during traction as well as the regenerative braking power decrease when the \gls{acr:em} temperature limit is reached. In general, the velocity profile shows smooth behavior (i.e., the acceleration of the vehicle is within a relatively small value range), which is typical for energy-limited scenarios~\cite{HerrmannChristEtAl2019,HerrmannPassigatoEtAl2020,LiuFotouhiEtAl2020}, since regenerative braking is used instead of the mechanical brakes to reduce the velocity before cornering. 
Third, we observe that the battery energy exceeds the terminal value before the end of the stint. Since the battery is partially charged in this scenario, the absolute lower energy limit is not reached. With the use of regenerative braking, the battery energy then reaches the terminal value exactly at the end of the stint, indicating an energy-limited scenario.
\revTCST{\begin{figure}[!tb]
	\centering
	\includegraphics[width=1\columnwidth,trim={10 20 35 41},clip]{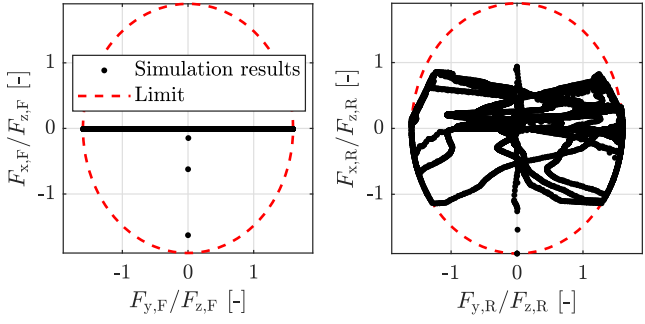}
	\caption{The friction circles per axle showing the normalized longitudinal and lateral forces. The vehicle is rear-wheel driven and braking is done mostly using the \gls{acr:em}. Therefore, the front axle (left figure) mostly generates lateral forces, whereas the rear axle (right figure) shows more combined forces. }
	\label{fig:friction_circles}
\end{figure}}
Finally, Fig.~\ref{fig:friction_circles} shows the normalized lateral and longitudinal forces per axle, along with the maximum grip limit defined by the friction coefficient. Since the vehicle in this study is rear-wheel driven \revTCST{($F_\mathrm{p,F}=0$)}, the front axle only provides negative forces in the longitudinal direction. It can be noticed that the data points do not show the typical diamond-like pattern for the front axle forces that is expected in racing scenarios~\cite{LovatoMassaro2021}. This is because the mechanical brakes are only used during pit entry at the end of the stint and because the vehicle applies regenerative braking to slow down for cornering. Moreover, both the front and rear axles are operated at the lateral grip limit in corners, indicating that the acceleration throughout the corners is maximized.

\subsection{High-level Optimization}
This section presents the optimal race strategy in terms of number of pit stops, stint length and charge time, and we compare it against the strategies optimized for a fixed number of pit stops.
We select a \unit[6]{h} race, yet longer races can be solved as well with our approach, considering the very low computational times needed by our high-level framework to converge.
To link the terminal battery energy $E_\mathrm{b}(S_\mathrm{stint})$ to the charge time $t_\mathrm{charge}$, we apply constant current charging, which is in line with the implementation on the actual vehicle. Furthermore, we calculate the terminal battery temperature $\vartheta_\mathrm{b,N}$ using the losses that correspond to the constant current profile and obtain the initial \gls{acr:em} temperature $\vartheta_\mathrm{m,init}$ by applying maximum cooling power, for the given charge time. 

\begin{figure}[!tb]
	\centering
	\input{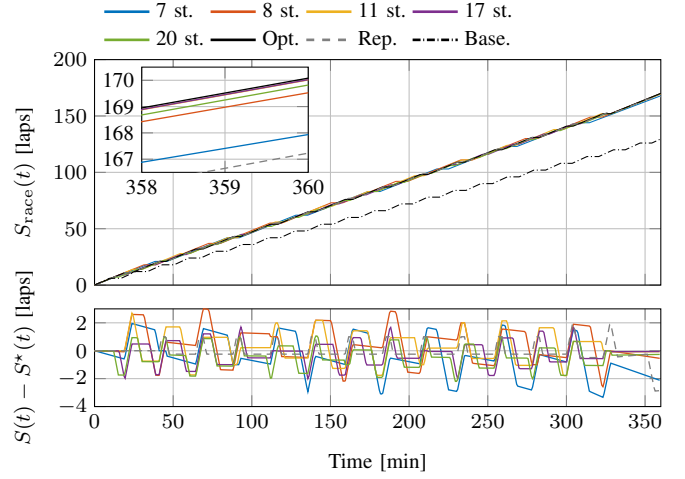}
	\caption{Evolution of the completed laps as a function of time for the optimal strategy (black) and the strategies optimized for a fixed number of pit stops. The dashed data show a repetition strategy of the global optimal stint and the dash-dotted data show a baseline strategy. The zoom window corresponds to the final 2 min of the race and illustrates the difference in race distance between the optimized strategies, showing that jointly optimizing the number of stints can significantly outperform other strategies by multiple laps.}
	\label{fig:cumulative_strategy}
\end{figure}

\par 
Fig.~\ref{fig:cumulative_strategy} shows the evolution of the completed laps as a function of time for various fixed pit stop strategies. We observe that the optimal strategy of \revTCST{\unit[14]{stops}} results in the largest amount of completed laps, thereby confirming that it is indeed optimal in terms of number of pit stops. The difference in covered race length between the optimal and fixed-pit-stop-number strategies can exceed multiple laps and hence significantly affect the final race outcome in terms of finishing position, highlighting the importance of jointly optimizing the number of pit stops. 
Furthermore, we compare the optimal race strategy against a strategy that repeats the global optimal stint of \revTCST{\unit[11]{laps} and \unit[4.7]{min} charging}, until the race ends. Although the stints used in this strategy are globally optimal from a stint perspective, it does not cover the largest distance within the \unit[6]{h} time window, which makes the strategy sub-optimal from a race perspective. This is because a pure repetition of the optimal stint does not perfectly fit within the \unit[6]{h} race. Instead, it is beneficial to slightly deviate from the optimal stint, so that the stints fit better within the race, as done by the optimal race strategy. This shows that the optimal race strategy is not necessarily the same as the optimal stint, thereby highlighting the importance of jointly optimizing the stint lengths, charge times and number of pit stops. 
Lastly, to show the importance of the bi-level approach, we compare the optimal strategy against a baseline \emph{flat-out} strategy whereby no energy management is applied to limit energy consumption, but the EMs are rather operated at maximum power whenever possible. This results in \unit[6]{lap} stints and a total race distance of \unit[132]{laps}, whilst the globally optimal solution is about \unit[170]{laps}. Hence, the optimal solution is about 29\% better compared to the baseline, even though the baseline was computed without thermal limits. 

\begin{figure}[!tb]
	\centering
	\input{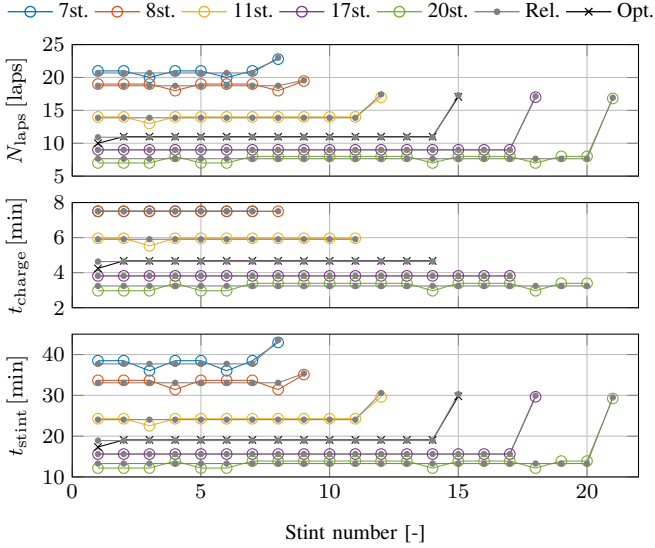}
	\caption{Optimal race strategy (black) in terms of stint length, charge time and stint time with $t_{\mathrm{charge,max}}=\unit[7.5]{min}$. For comparison, we show other optimal fixed pit stops strategies together with the relaxed solution in gray. Stint length, charge- and stint time are related and the optimal integer solution minimizes the differences to the relaxed solution.}
	\label{fig:stint_strategy}
\end{figure}

Fig.~\ref{fig:stint_strategy} shows the individual stints in terms of length and charge time, together with the relaxed non-integer solution.
We can conclude that a constant stint length over the race is optimal, since all stints in the relaxed solutions are equal, the only exception being the last stints.
In this use case, the optimal integer solution consists of the stint lengths that minimize the difference to the relaxed solution, namely, of a stint length \revTCST{between 10 and \unit[11]{laps} together with a charge time of about \unit[4.7]{min} and 14 pit stops} in total.
For the vehicle considered in this study, the \unit[8]{stop} strategy is the fastest strategy that involves fully charging the battery. However, this strategy is \unit[0.6]{laps} behind on the optimal race strategy, thereby showing that fully charging the battery can be sub-optimal. This is explained by the constant current charging, which effectively reduces the average charging power for longer charge times due to the battery voltage decreasing for lower energy levels. With a different charging profile, such as constant power charging, a different optimal race strategy will be obtained. However, it is beyond the scope of this paper to optimize the charging profile.
Moreover, reducing the charge time results in shorter stints, thereby increasing the number of pit stops needed, as indicated by the results. This then results in more time lost in the pit lane due to the speed limit. Thus there is a clear trade-off in determining the optimal charge time. When we also consider thermally-limited scenarios, this trade-off becomes even more complex. In these cases, a shorter charge time is expected to be in favor, since this increases the terminal battery temperature. Although a longer charge time reduces the initial \gls{acr:em} temperature, we observed that this had a relatively small impact on the stint time. In fact, it has been observed that in longer stints the \gls{acr:em} was operated relatively more at the thermal limit compared to shorter stints, thereby disfavoring longer charge times. Lastly, we observe that there is a considerable decrease in race distance when an increasing stint length cannot be compensated by an increase in charge time, as illustrated by the difference between the 7 and 8 stop strategies. 
From the aforementioned observations, we conclude that the stint length, stint time and charge time are closely related in the case of an optimal solution.

\subsection{Validation}
In this section, we validate the correctness of the model by showing that the lethargy constraint holds with equality for various stints that were presented in the previous section. Furthermore, we validate the convex models by implementing the optimal inputs into a non-linear simulator and compare the drift in battery energy. Finally, we validate the numerical combinations of stint length and charge time for the various strategies. For the latter, we calculate the average stint velocity $\overline{v}_{\mathrm{stint}}(k)$ for every strategy as
\par\nobreak\vspace{-5pt}
\begingroup
\allowdisplaybreaks
\begin{small}
	\begin{equation}\label{eq:stint_velocity}
		\overline{v}_{\mathrm{stint}}(k) = \frac{S_{\mathrm{stint}}(k)}{t_{\mathrm{charge}}(k)+t_{\mathrm{stint}}(k)}, \quad \forall k \in [0,n_\mathrm{stops}-1].
	\end{equation}
\end{small}%
\endgroup
\revTCST{\begin{figure}[!tb]
	\vspace{2mm}
	\centering
	\small
	\includegraphics[width=1\columnwidth,trim={10 0 15 5},clip]{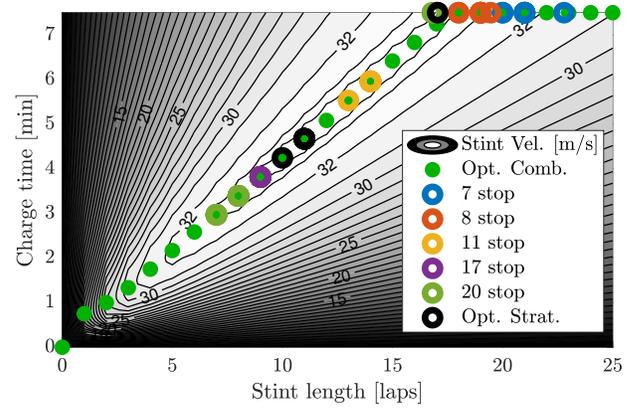}
	\caption{The average stint velocity for a combination of stint lengths and charge times together with the optimal combinations and actual numerical solutions. The optimal combinations of stint length and charge time show a clear (linear) correlation, to which the numerical solutions are aligned.}
	\label{fig:stint_velocity}
\end{figure}}
The globally optimal stint should maximize the average stint velocity, since that maximizes the driven distance per unit of time, \revTCST{ also accounting for the charge time after the stint}. \revTCST{We can calculate the average stint velocity for various combinations of stint length and charge time using the stint time map.} \revTCST{To validate the correctness of the optimized strategies, we  calculate the optimal charge time that maximizes the average stint velocity, for every possible stint length. These are referred to as the optimal combinations of charge times and stint length. The optimized strategies should align with these theoretical optimal combinations when looking from a stint perspective. However, small deviations might be observed, since the strategies are optimized from a race perspective. Fig.~\ref{fig:stint_velocity} shows the average stint velocity together with the theoretical optimal combinations.} These optimal combinations indicate an almost linear relation between charge time and stint length until the maximum charge time is reached. The most noticeable exception is the single-lap stint, which has a relatively higher optimal charge time. This is because a single-lap stint does not include a flying lap, where the vehicle starts the lap with a high initial velocity. Thus, the vehicle requires more energy in the first lap to accelerate towards the optimal velocity, which explains the longer charge time for a single-lap stint.
The globally optimal stint \revTCST{with the highest average stint velocity} consists of \revTCST{ \unit[11]{laps} and \unit[4.7]{min} charging}, which is the exact same combination that we obtained as the optimal strategy in the previous section. Furthermore, we observe that the average stint velocity decreases in sensitivity around the optimal combinations for increasing stint length and charge time. When the maximum charge time is reached, the average stint velocity diminishes considerably for increasing stint lengths. Thereby, increasing the stint length beyond \unit[18]{laps} quickly becomes less favorable, since it cannot be compensated by an increase in charge time. This explains why the \unit[7]{stop} strategy is significantly worse than the others. 
Finally, we note that the numerical solutions are in line with the theoretically optimal combinations, \revTCST{thereby validating their correctness}.
The outliers not aligning with the optimal combinations, e.g., at \revTCST{ \unit[16.8]{laps} and \unit[7.5]{min} charging}, correspond to the last stints, for which the charge time is not part of the race and thus the calculation of the stint velocity in~\eqref{eq:stint_velocity} is not valid.

\revTCST{\begin{figure}[!tb]
	\vspace{2mm}
	\centering
	\small
	\revTCST{\includegraphics[width=\columnwidth,trim={10 0 18 0},clip]{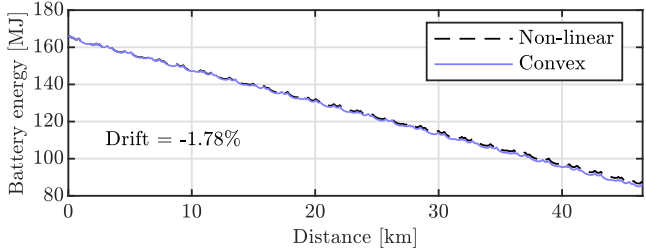}}
	\caption{Battery energy calculated through the non-linear and convex models, using the optimized inputs for the globally optimal stint of \revTCST{\unit[11]{laps} and \unit[4.7]{min} charging}. The small drift in battery energy between the non-linear and convex models indicate that the convex models accurately capture the dynamics of the powertrain components.}
	\label{fig:NLsimulator}
\end{figure}}

To validate the convex loss models for the \gls{acr:em} and battery, we calculate the battery energy trajectory by applying the optimal input trajectories to non-linear models and compare the result to the trajectory obtained from the optimization. Fig.~\ref{fig:NLsimulator} shows both battery energy trajectories for the globally optimal stint of \revTCST{\unit[11]{laps} with \unit[4.7]{min}} charging. From this figure, we observe small deviations, with a total drift of \revTCST{ \unit[-1.78]{\%}} with respect to the non-linear models, thereby indicating that the convex models accurately capture the dynamics of the powertrain components.

\revTCST{\begin{figure}[!tb]
	\vspace{2mm}
	\centering
	\small
	\includegraphics[width=\columnwidth]{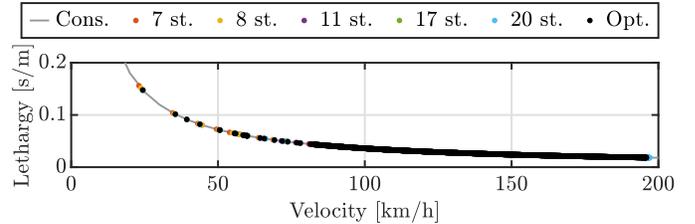}
	\caption{Validation of the lethargy constraint for the mode of the stints obtained from the high-level optimization. All data points align with the constraint (within the solver tolerance), indicating that the lethargy constraint holds with equality.}
	\label{fig:val_lethargy}
\end{figure}}

Finally, we verify that all relaxed constraints hold with equality. As noted previously, it can occur that the lethargy constraint does not hold with equality in thermally-limited scenarios. Therefore, we explicitly show this constraint for the most frequently used stints per pit stop strategy that were obtained in the previous section. Fig.~\ref{fig:val_lethargy} shows that all data points align with the constraint, thereby indicating that the lethargy constraint holds with equality for all optimal stints. 

\section{Conclusion} \label{Conclusion}
In this paper, we studied the maximum-distance endurance racing problem for electric race cars, explicitly accounting for temperature and vehicle dynamics. To this end, we devised a convex bi-level optimization framework that can efficiently compute the optimal race strategy (in terms of stint lengths, charge times and pit stops), jointly with the in-stint \gls{acr:PO} strategy.  

Our results showed that accounting for the temperature dynamics of the powertrain components is extremely important in an endurance racing scenario and that both the in-stint \gls{acr:PO} and the high-level race strategy were affected by it. 
From a stint perspective, there is a clear correlation between optimal stint length and charge time, which corresponds to the maximization of the average stint velocity. 
In contrast to our preliminary results~\cite{KampenHerrmannEtAl2022}, fully charging the battery proved to be sub-optimal for this case study, because increasing the battery operating window results in a lower average charge power and stricter battery temperature limitations. 
Moreover, the results suggest to avoid the usage of the mechanical friction brakes, which indicates the paramount importance of efficient energy management for electric race cars. In fact, this highlights the shortcomings and need for improvement of current battery technologies, compared to hybrid-electric vehicles in racing scenarios.   
Finally, we highlighted the importance of optimizing both levels by showing that a repetition of the optimal stint does not result in an optimal race.

This work opens the field for the following possible extensions: 
First, we are interested in the impact of tire degradation \revTCST{and temperature} on the achievable stint time and the resulting race strategies, \revTCST{ since this yields an apparent trade-off with energy management. }
Second, we want to extend the usage of the framework to an investigation of the optimal transmission technology and powertrain architecture.
Finally, we would like to leverage the computational efficiency of the framework for an online implementation both in the vehicle, as well as in the pit lane.

\ifextendedversion
\section*{Appendix} \label{sec:Appendix}
In Section~\ref{sec:Methodology - low-level optimization}, we frequently used the following relaxations
\par\nobreak\vspace{-5pt}
\begingroup
\allowdisplaybreaks
\begin{small}
	\begin{equation}
	x\cdot y\geq z^\top z,
		\label{eq:appendix1}
	\end{equation}
\end{small}%
\endgroup
which defines a convex set. To prove convexity, we can rewrite this constraint to a convex second-order conic constraint as
\par\nobreak\vspace{-5pt}
\begingroup
\allowdisplaybreaks
\begin{small}
	\begin{equation}
		x + y\geq 
		\begin{Vmatrix}
			2\cdot z \\
		x-y\\	
		\end{Vmatrix}_2,
		\label{eq:appendix_socc}
	\end{equation}
\end{small}%
\endgroup
which can be solved with global optimality guarantees~\cite{BoydVandenberghe2004}. Since~\eqref{eq:appendix1} is mathematically equivalent to~\eqref{eq:appendix_socc}, both optimization problems will converge to the same KKT points, thereby guaranteeing global optimality. 

Another type of constraint that was used in Section~\ref{sec:Methodology - low-level optimization} for the friction circles is
\par\nobreak\vspace{-5pt}
\begingroup
\allowdisplaybreaks
\begin{small}
	\begin{equation}
		x^2+ y^2 \leq z^2.
		\label{eq:appendix2}
	\end{equation}
\end{small}%
\endgroup
This constraint can be directly translated to a convex second-order conic constraint
\par\nobreak\vspace{-5pt}
\begingroup
\allowdisplaybreaks
\begin{small}
	\begin{equation}
			\begin{Vmatrix}
			x \\
			y\\	
		\end{Vmatrix}_2 \leq z,
		\label{eq:appendix6}
	\end{equation}
\end{small}%
\endgroup
thereby proving that~\eqref{eq:appendix2} defines a convex set.

In Section~\ref{sec:results - low-level} we mentioned that an active terminal battery temperature constraint results in a gradual decrease in battery temperature over the entire stint. Fig.~\ref{fig:battery_limited} again shows an \unit[11]{lap} stint, but with an increased thermal resistance for the battery, thereby resulting in an active battery temperature constraint. Comparing this figure to Fig.~\ref{fig:Low_level_results}, we can clearly see the difference in battery temperature trajectories. Whereas the battery temperature was kept close to the maximum before to maximize the efficiency, the temperature is now gradually decreasing to reach the required terminal value.

\begin{figure}[!tb]
	\centering
	\input{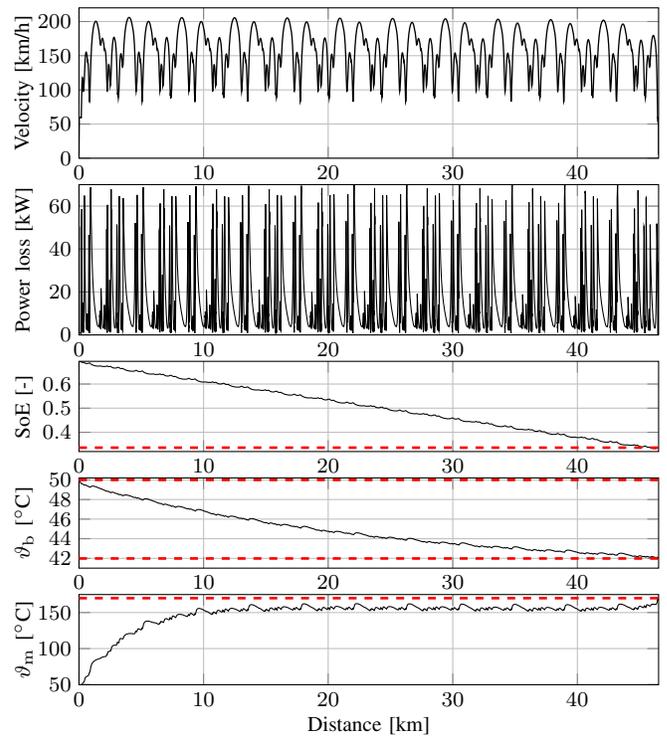}
	\caption{Velocity, powertrain losses, battery SoE, battery temperature and \gls{acr:em} temperature trajectories for an 11 lap stint. The terminal battery temperature is an active constraint, which results in a gradual decrease in battery temperature throughout the stint to allow for some margin in temperature for the consecutive charging session.}
	\label{fig:battery_limited}
\end{figure}
\else

\fi

\section*{Acknowledgment}
\noindent We thank Dr.~I.~New for proofreading this paper, and N.~Schoenmakers for the photograph of the car shown in Fig.~\ref{fig:inmotion-car}. This paper was partly supported by the NEON research project (project number 17628 of the Crossover program which is (partly) financed by the Dutch Research Council (NWO)).

\bibliographystyle{IEEEtran}
\bibliography{main,SML_papers}

\newcommand{\noopsort}[1]{} \newcommand{\printfirst}[2]{#1}
  \newcommand{\singleletter}[1]{#1} \newcommand{\switchargs}[2]{#2#1}
\begin{thebibliography}{10}
\providecommand{\url}[1]{#1}
\csname url@samestyle\endcsname
\providecommand{\newblock}{\relax}
\providecommand{\bibinfo}[2]{#2}
\providecommand{\BIBentrySTDinterwordspacing}{\spaceskip=0pt\relax}
\providecommand{\BIBentryALTinterwordstretchfactor}{4}
\providecommand{\BIBentryALTinterwordspacing}{\spaceskip=\fontdimen2\font plus
\BIBentryALTinterwordstretchfactor\fontdimen3\font minus
  \fontdimen4\font\relax}
\providecommand{\BIBforeignlanguage}[2]{{%
\expandafter\ifx\csname l@#1\endcsname\relax
\typeout{** WARNING: IEEEtran.bst: No hyphenation pattern has been}%
\typeout{** loaded for the language `#1'. Using the pattern for}%
\typeout{** the default language instead.}%
\else
\language=\csname l@#1\endcsname
\fi
#2}}
\providecommand{\BIBdecl}{\relax}
\BIBdecl

\bibitem{InMotion}
(2023) {InMotion fully electric LMP3 car}. {InMotion}. Available at
  \url{https://www.inmotion.tue.nl/en/about-us/cars/revolution}.

\bibitem{FIA2021}
FIA. (2021) 2021 fia world endurance championship: Sporting regulations.
  Available online at
  \url{https://www.fia.com/sites/default/files/2021_wec_sporting_regulations_-_wmsc161220-clean_1.pdf}.

\bibitem{LotEvangelou2013}
R.~Lot and S.~Evangelou, ``Lap time optimization of a sports series hybrid
  electric vehicle,'' in \emph{{World Congress on Engineering}}, 2013.

\bibitem{SedlacekOdenthalEtAl2020}
T.~Sedlacek, D.~Odenthal, and D.~Wollherr, ``Minimum-time optimal control for
  battery electric vehicles with four wheel-independent drives considering
  electrical overloading,'' \emph{{Vehicle System Dynamics}}, 2020.

\bibitem{Casanova2000}
D.~Casanova, ``On minimum time vehicle manoeuvring: The theoretical optimal
  lap,'' Ph.D. dissertation, {Cranfield University}, 2000.

\bibitem{LimebeerPerantoni2014}
D.~Limebeer and G.~Perantoni, ``Optimal control for a formula one car with
  variable parameters,'' \emph{{Vehicle System Dynamics}}, vol.~52, no.~5, pp.
  653--678, 2014.

\bibitem{ChristWischnewskiEtAl2019}
F.~Christ, A.~Wischnewski, A.~Heilmeier, and B.~Lohmann, ``Time-optimal
  trajectory planning for a race car considering variable tyre-road friction
  coefficients,'' \emph{{Vehicle System Dynamics}}, vol.~59, pp. 588--612,
  2019.

\bibitem{DalBiancoLotEtAl2017}
N.~Dal~Bianco, R.~Lot, and M.~Gadola, ``Minimum time optimal control simulation
  of a gp2 race car,'' \emph{{Journal of Automobile Engineering}}, vol. 232,
  pp. 1180--1195, 2017.

\bibitem{HeilmeierWischnewskiEtAl2020}
A.~Heilmeier, A.~Wischnewski, L.~Hermansdorfer, J.~Betz, M.~Lienkamp, and
  B.~Lohmann, ``Minimum curvature trajectory planning and control for an
  autonomous race car,'' \emph{{Vehicle System Dynamics}}, vol.~58, no.~10, pp.
  1497--1527, 2020.

\bibitem{MassaroLimebeer2021}
M.~Massaro and D.~J.~N. Limebeer, ``Minimum-lap-time optimisation and
  simulation,'' \emph{{Vehicle System Dynamics}}, 2021.

\bibitem{TremlettLimebeer2016}
A.~J. Tremlett and D.~J.~N. Limebeer, ``Optimal tyre usage for a formula one
  car,'' \emph{{Vehicle System Dynamics}}, 2016.

\bibitem{HerrmannChristEtAl2019}
T.~Herrmann, F.~Christ, J.~Betz, and M.~Lienkamp, ``Energy management strategy
  for an autonomous electric racecar using optimal control,'' in \emph{{Proc.\
  IEEE Int.\ Conf.\ on Intelligent Transportation Systems}}, 2019.

\bibitem{LimebeerPerantoniEtAl2014}
D.~Limebeer, G.~Perantoni, and A.~Rao, ``Optimal control of formula one car
  energy recovery systems,'' \emph{{Int.\ Journal of Control}}, vol.~87,
  no.~10, pp. 2065--2080, 2014.

\bibitem{HerrmannPassigatoEtAl2020}
T.~Herrmann, F.~Passigato, J.~Betz, and M.~Lienkamp, ``Minimum race-time
  planning-strategy for an autonomous electric racecar,'' in \emph{{Proc.\ IEEE
  Int.\ Conf.\ on Intelligent Transportation Systems}}, 2020.

\bibitem{LiuFotouhiEtAl2020}
X.~Liu, A.~Fotouhi, and D.~Auger, ``Optimal energy management for formula-e
  cars with regulatory limits and thermal constraints,'' \emph{{Applied
  Energy}}, vol. 279, 2020.

\bibitem{HerrmannSauerbeckEtAl2021}
T.~Herrmann, F.~Sauerbeck, M.~Bayerlein, and M.~Betz, J. an~Lienkamp,
  ``Optimization-based real-time-capable energy strategy for autonomous
  electric race cars,'' \emph{{SAE International Journal of Connected and
  Automated Vehicles}}, 2021, in press.

\bibitem{LiuFotouhi2020}
X.~Liu and A.~Fotouhi, ``Formula-e race strategy development using artificial
  neural networks and monte carlo tree search,'' \emph{{Neural Computing and
  Applications}}, no.~32, p. 15191–15207, 2020.

\bibitem{EbbesenSalazarEtAl2018}
S.~Ebbesen, M.~Salazar, P.~Elbert, C.~Bussi, and C.~H. Onder, ``Time-optimal
  control strategies for a hybrid electric race car,'' \emph{{IEEE Transactions
  on Control Systems Technology}}, vol.~26, no.~1, pp. 233--247, 2018.

\bibitem{SalazarElbertEtAl2017}
M.~Salazar, P.~Elbert, S.~Ebbesen, C.~Bussi, and C.~H. Onder, ``Time-optimal
  control policy for a hybrid electric race car,'' \emph{{IEEE Transactions on
  Control Systems Technology}}, vol.~25, no.~6, pp. 1921--1934, 2017.

\bibitem{DuhrChristodoulouEtAl2020}
P.~Duhr, G.~Christodoulou, C.~Balerna, M.~Salazar, A.~Cerofolini, and C.~H.
  Onder, ``Time-optimal gearshift and energy management strategies for a hybrid
  electric race car,'' \emph{{Applied Energy}}, vol. 282, no. 115980, 2020.

\bibitem{BorsboomFahdzyanaEtAl2021}
O.~Borsboom, C.~A. Fahdzyana, T.~Hofman, and M.~Salazar, ``A convex
  optimization framework for minimum lap time design and control of electric
  race cars,'' \emph{{IEEE Transactions on Vehicular Technology}}, vol.~70,
  no.~9, pp. 8478--8489, 2021.

\bibitem{LocatelloKondaEtAl2020}
A.~Locatello, M.~Konda, O.~Borsboom, T.~Hofman, and M.~Salazar, ``Time-optimal
  control of electric race cars under thermal constraints,'' in \emph{{European
  Control Conference}}, 2021.

\bibitem{HeilmeierGrafEtAl2018}
A.~Heilmeier, M.~Graf, and M.~Lienkamp, ``A race simulation for strategy
  decisions and circuit motorsports,'' in \emph{{Proc.\ IEEE Int.\ Conf.\ on
  Intelligent Transportation Systems}}, 2018.

\bibitem{WestLimebeer2020}
W.~J. West and D.~J.~N. Limebeer, ``Optimal tyre management for a
  high-performance race car,'' \emph{{Vehicle System Dynamics}}, vol. 231, pp.
  1--19, 2020.

\bibitem{KampenHerrmannEtAl2022}
J.~van Kampen, T.~Herrmann, and M.~Salazar, ``Maximum-distance race strategies
  for a fully electric endurance race car,'' \emph{{European Journal of
  Control}}, 2022, in press. Available online at
  \url{https://arxiv.org/abs/2111.05784}.

\bibitem{BroereSalazar2022}
S.~Broere, J.~van Kampen, and M.~Salazar, ``Minimum-lap-time control strategies
  for all-wheel drive electric race cars via convex optimization,'' in
  \emph{{European Control Conference}}, 2022.

\bibitem{GuzzellaSciarretta2007}
L.~Guzzella and A.~Sciarretta, \emph{Vehicle propulsion systems: Introduction
  to Modeling and Optimization}, 2nd~ed.\hskip 1em plus 0.5em minus 0.4em\relax
  {Springer Berlin Heidelberg}, 2007.

\bibitem{BorsboomFahdzyanaEtAl2020}
O.~Borsboom, C.~A. Fahdzyana, M.~Salazar, and T.~Hofman, ``Time-optimal control
  strategies for electric race cars for different transmission technologies,''
  in \emph{{IEEE Vehicle Power and Propulsion Conference}}, 2020.

\bibitem{ZhangXiaEtAl2018}
R.~Zhang, B.~L. Bizhong Xia~and, L.~Cao, Y.~Lai, W.~Zheng, H.~Wang, W.~Wang,
  and M.~Wang, ``A study on the open circuit voltage and state of charge
  characterization of high capacity lithium-ion battery under different
  temperature,'' \emph{Energies}, 2018.

\bibitem{FarmannSauer2017}
A.~Farmann and D.~U. Sauer, ``A study on the dependency of the open-circuit
  voltage on temperature and actual aging state of lithium-ion batteries,''
  \emph{Journal of Power Sources}, vol. 347, 2017.

\bibitem{RosewaterCoppEtAl2019}
D.~M. Rosewater, D.~A. Copp, T.~A. Nguyen, R.~H. Byrne, and S.~Santoso,
  ``Battery energy storage models for optimal control,'' \emph{IEEE Access},
  vol.~7, pp. 178\,357--178\,391, 2019.

\bibitem{Lebkowski2017}
A.~Łebkowski, ``Temperature, overcharge and short-circuit studies of batteries
  used in electric vehicles,'' \emph{Przeglad Elektrotechniczny}, 2017.

\bibitem{BoydVandenberghe2004}
S.~Boyd and L.~Vandenberghe, \emph{Convex optimization}.\hskip 1em plus 0.5em
  minus 0.4em\relax {Cambridge Univ.\ Press}, 2004.

\bibitem{RichardsHow2005}
A.~Richards and J.~How, ``Mixed-integer programming for control,'' in
  \emph{{Proc.\ of the American Control Conference}}, 2005.

\bibitem{Lee2012}
J.~Lee and S.~Leyffer, Eds., \emph{Mixed Integer Nonlinear Programming}.\hskip
  1em plus 0.5em minus 0.4em\relax {Springer-Verlag}, 2012.

\bibitem{BelottiKirchesEtAl2013}
P.~Belotti, C.~Kirches, S.~Leyffer, J.~Linderoth, J.~Luedtke, and A.~Mahajan,
  ``Mixed-integer nonlinear optimization,'' \emph{Acta Numerica}, 2013.

\bibitem{AnderssonGillisEtAl2019}
J.~A.~E. Andersson, J.~Gillis, G.~Horn, J.~B. Rawlings, and M.~Diehl, ``Casadi
  -- a software framework for nonlinear optimization and optimal control,''
  \emph{{Mathematical Programming Computation}}, vol.~11, no.~1, pp. 1--36,
  2019.

\bibitem{WachterBiegler2006}
A.~Wachter and L.~T. Biegler, ``On the implementation of an interior-point
  filter line-search algorithm for large-scale nonlinear programming,''
  \emph{{Mathematical Programming}}, vol. 106, no.~1, pp. 25--57, 2006.

\bibitem{HSL}
\BIBentryALTinterwordspacing
HSL. A collection of fortran codes for large scale scientific computation.
  [Online]. Available: \url{http://www.hsl.rl.ac.uk}
\BIBentrySTDinterwordspacing

\bibitem{Loefberg2004}
J.~L{\"o}fberg, ``{YALMIP} : A toolbox for modeling and optimization in
  {MATLAB},'' in \emph{{IEEE Int.\ Symp.\ on Computer Aided Control Systems
  Design}}, 2004.

\bibitem{ApS2017}
M.~{ApS}. (2017) {MOSEK} optimization software. {Available at
  }\url{https://mosek.com/}.

\bibitem{LovatoMassaro2021}
S.~Lovato and M.~Massaro, ``A three-dimensional free-trajectory
  quasi-steady-state optimal-control method for minimum-lap-time of race
  vehicles,'' \emph{{Vehicle System Dynamics}}, 2021.

\end{thebibliography}

\begin{IEEEbiography}[{\includegraphics[width=1in,height=1.25in,clip,keepaspectratio]{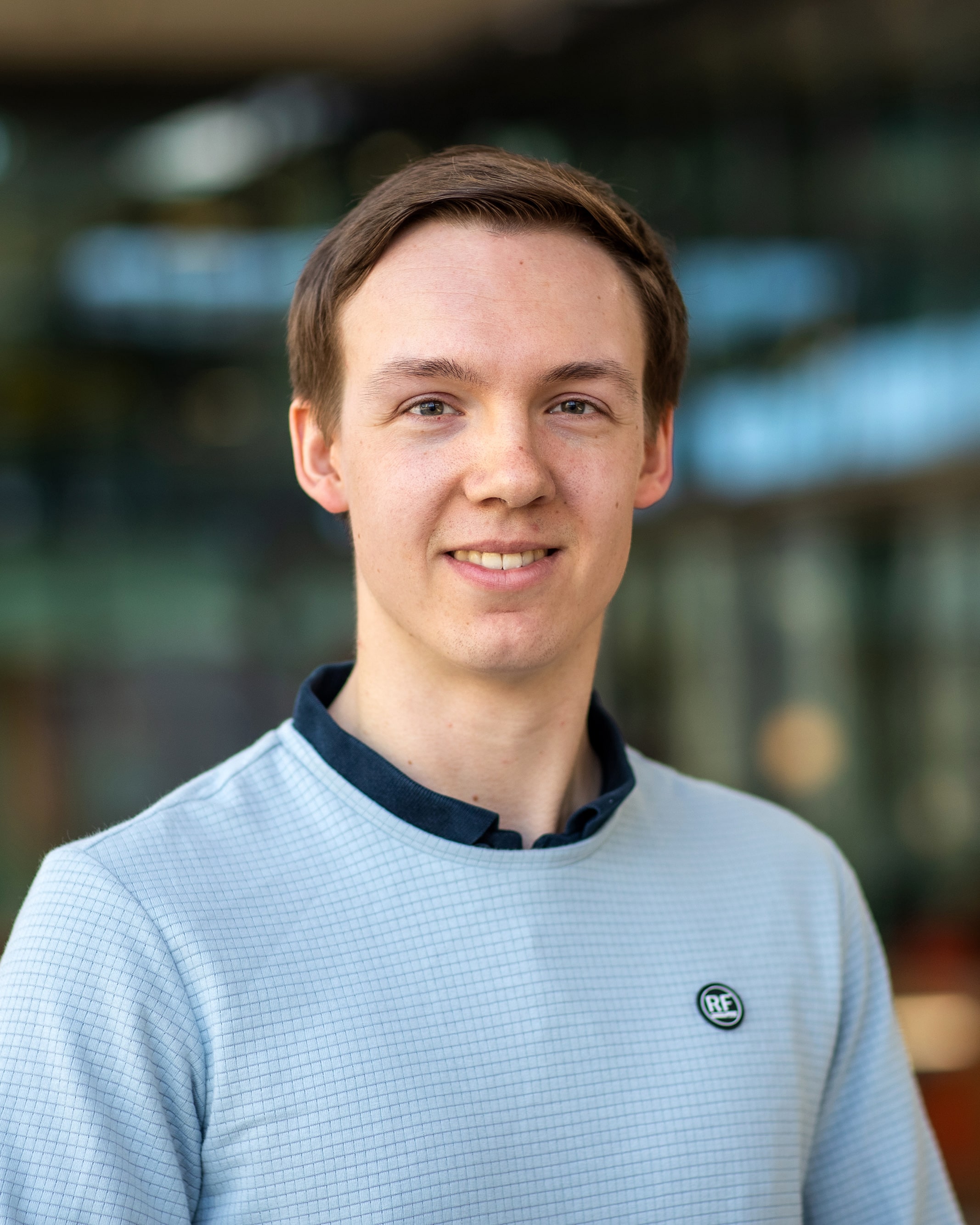}}]{Jorn van Kampen} obtained the B.Sc. degree (with great appreciation) and the M.Sc. degree (cum laude) in Automotive Technology from Eindhoven University of Technology (TU/e), The Netherlands, in 2018 and 2022, respectively. He is currently pursuing a Ph.D. degree at the Control Systems Technology section of TU/e.
His main research interests include electric race vehicles and optimization methods for powertrain design and control.
He was granted the Best Student Paper Award at the 2022 European Control Conference and his Master thesis received the 2022 Tata Steel Graduation Award for Mechanical Engineering and Material Science.
\end{IEEEbiography}

\begin{IEEEbiography}[{\includegraphics[width=1in,height=1.25in,clip,keepaspectratio,trim = 0 150 0 20]{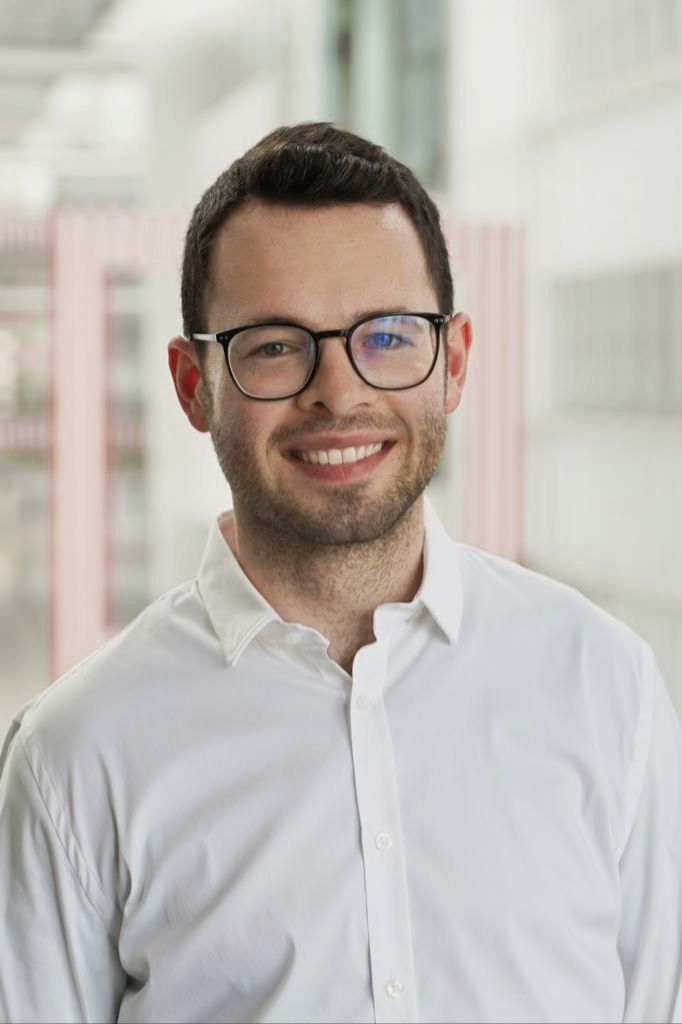}}]{Thomas Herrmann} was awarded a B.Sc. and a	M.Sc. in Mechanical Engineering by the Technical
	University of Munich (TUM), Germany, in 2016
	and 2018, respectively. From 2018 – 2022 he was pursuing
	his doctoral studies at the Institute of Automotive
	Technology at TUM. His research interests
	include optimal control in the field of trajectory
	planning for autonomous vehicles and the efficient
	incorporation of the electric powertrain behavior
	within these optimization problems. 
	He was granted the Best Student Paper Award at the 2022 European Control Conference.
\end{IEEEbiography}

\begin{IEEEbiography}[{\includegraphics[width=1in,height=1.25in,clip,keepaspectratio,trim = 70 100 70 30]{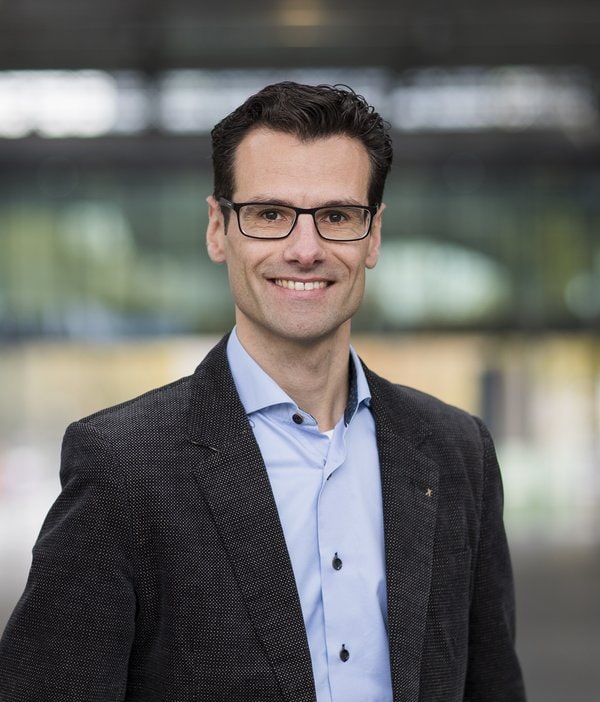}}]{Theo Hofman} was born in Utrecht, The Netherlands, in 1976. He received the M.Sc. (with Hons.) and Ph.D. degrees in mechanical engineering from the Eindhoven University of Technology, Eindhoven, The Netherlands, in 1999 and 2007, respectively. From 1999 to 2003, he was a Researcher and Project Manager with the R\&D Department of Thales–Cryogenics B.V., Eindhoven, The Netherlands. From 2003 to 2007, he was a Scientific Researcher with Drivetrain Innovations B.V., Eindhoven, The Netherlands. Since 2010, he has been an Associate Professor with the Control Systems Technology Group. His research interests include system design optimization methods for complex dynamical engineering systems and discrete topology design using computational design synthesis.
\end{IEEEbiography}

\begin{IEEEbiography}[{\includegraphics[width=1in,height=1.25in,clip,keepaspectratio]{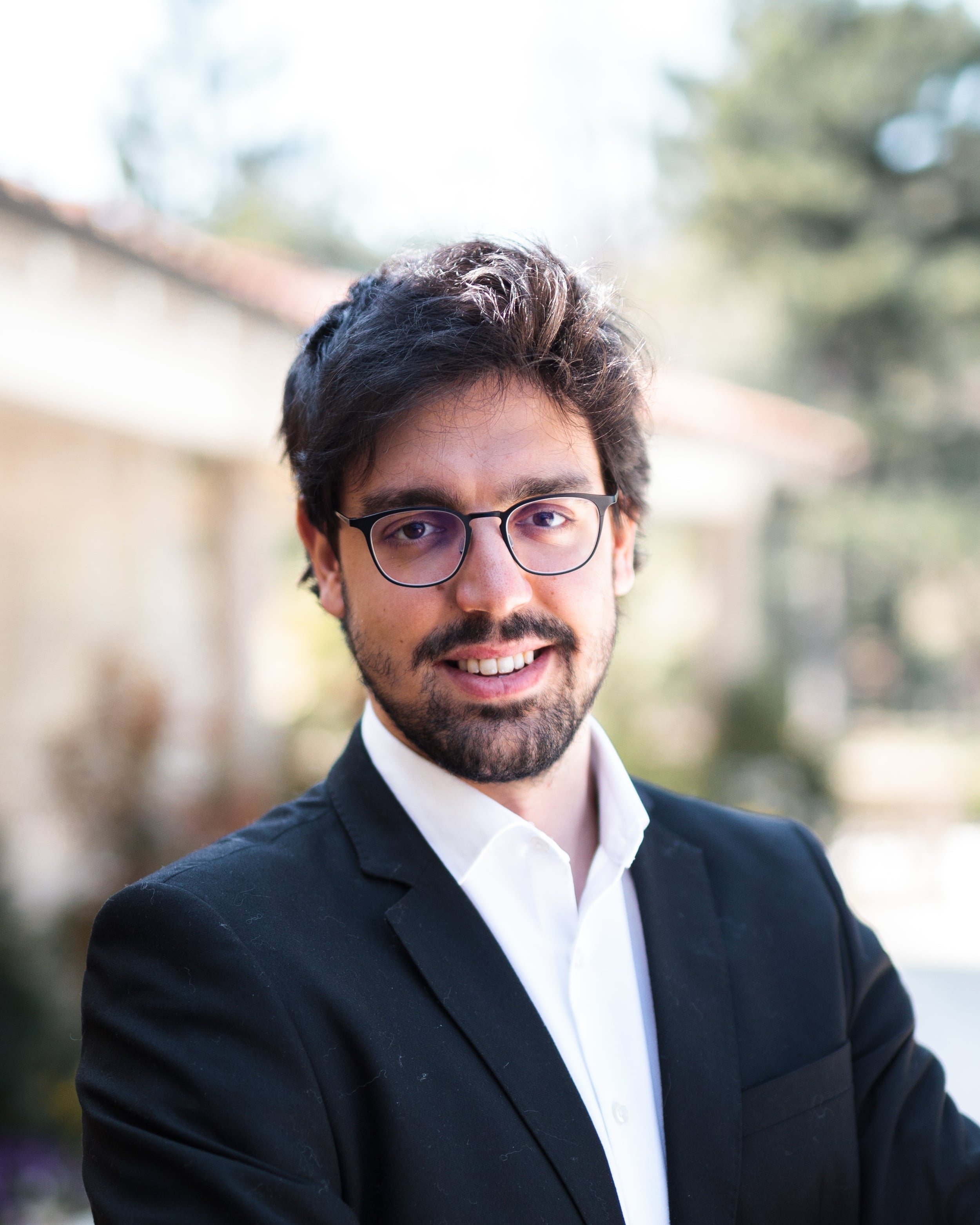}}]{Mauro Salazar} is an Assistant Professor in the Control Systems Technology section at Eindhoven University of Technology (TU/e), and co-affiliated with Eindhoven AI Systems Institute (EAISI). He received the Ph.D. degree in Mechanical Engineering from ETH Zurich in 2019. Before joining TU/e he was a Postdoctoral Scholar in the Autonomous Systems Lab at Stanford University.
Dr. Salazar’s research is focused on optimization models and methods for cyber-socio-technical systems design and control, with a strong focus on sustainable mobility.
Both his Master thesis and PhD thesis were recognized with the ETH Medal, and his papers were granted the Best Student Paper award at the 2018 Intelligent Transportation Systems Conference and at the 2022 European Control Conference.
\end{IEEEbiography}

\end{document}